\documentclass[12pt]{article}
\date{}
\makeatletter
\@addtoreset{equation}{section}
\makeatother
\makeatother
\usepackage{amssymb,amsmath}
\usepackage{amsmath}
\usepackage{color}
\usepackage[colorlinks,linkcolor=blue,citecolor=blue]{hyperref}
\usepackage{cite}
\usepackage{graphicx}
\usepackage[misc]{ifsym}

\numberwithin{equation}{section}
\newtheorem{thm}{Theorem}[section]
\newtheorem{lem}[thm]{Lemma}
\newtheorem{rem}[thm]{Remark}

\newtheorem{cor}[thm]{Corollary}

\allowdisplaybreaks

\begin{document}
	\title{\Large Ground state solutions for the asymptotically	periodic Schr\"odinger-Poisson  systems with $p$-Laplacian	\footnote{	Supported by NSFC12301144 and 2024NSFSC1342.} }
	\author{\small  Yao Du\ \ \ \  Linfeng Fan\footnote{Corresponding author. \hfill\break\indent \ \ \ E-mail addresses:\ \ 1052591976@qq.com(Y. Du),  2090475403@qq.com(L. Fan).}
	\ \ \ \  \\
		{School of  Science, XiHua University}\\  {\small Chengdu 610039, People's Republic of China}}{\small }
		{\small } {\small }

	\maketitle
	\begin{abstract}
	In this paper we study the existence of ground state solutions for the asymptotically periodic Schr\"odinger-Poisson systems which are coupled by a Schr\"odinger equation of $p$-Laplacian and a Poisson equation of $q$-Laplacian. The method relies on a variational approach and the case of the nonlinearity exhibits	a critical growth is also considered. Some results in the literature are extended.
	\end{abstract}
	\section{Introduction}
~~~~This paper deal with the existence of ground state solutions for the quasilinear	Schr\"odinger-Poisson system
\begin{eqnarray} \label{1001}
	\left\{
	\begin{array}{ll}
		-\Delta_pu+V(x)|u|^{p-2}u+K(x)\phi|u|^{m-2}u=g(x,u)&\quad\text{in}\:\mathbb{R}^3,\\-\Delta_q\phi=K(x)|u|^{m}&\quad\text{in}\:\mathbb{R}^3,
	\end{array}
	\right.
\end{eqnarray}
where $V,K$ and $g$ are asymptotically periodic functions in $x$, $\Delta_iu=\text{div}(|\nabla u|^{i-2}\nabla u)$$\left(i=p,q\right)$, $1<p<3$, $q$ and $m$ satisfy
\begin{eqnarray}\label{102}
	\begin{array}{ll}\frac{3p}{4p-3}<q<3,\quad\frac{(4q-3)p}{3q}<m<\frac{(q-1)p^*}{q},	\end{array}
\end{eqnarray}
$p^*=\frac{3p}{3-p}$.
With $p=q=m=2$. System \eqref{1001} reduces to the semilinear Schr\"odinger-Poisson system
\begin{eqnarray} \label{0001500}
	\left\{
	\begin{array}{ll}-\Delta u + V(x)u + K(x)\phi u = g(x,u)&\quad \text{in~} \mathbb{R}^3,\\
		-\Delta \phi = K(x)u^2 &\quad \text{in~} \mathbb{R}^3,	\end{array}
	\right.
\end{eqnarray}
which has been studied intensively in the last twenty years due to its physical background, starting with the pioneering paper \cite{BF1998}.
For the results on existence and non-existence of solutions, radial and non-radial solutions, positive and sign-changing solutions, multiplicity of solutions, and ground states of system \eqref{0001500}, we refer to \cite{2011Alve,AAR2,Az2010,2025Murcia,AAP2008,GG2010,MA2003,DDM2004,2002ANSdA,DTM2004,DR2006,2016SunJT,WXG2017,WZGcritical2019,2014LiuZ,ZH2007,2008ZhaoLG,ZT2018,ZXZ2013,ZXZ2014,ZZ2009}.

Azzollini and   Pomponio \cite{AAP2008}  first   considered  the existence of a ground state solution for system \eqref{0001500} with $K=1$ and $g(x,u)=|u|^{r-1} u,$
when $V$ is a positive constant and $2 < r < 5$, or $V$ is possibly unbounded below satisfying $V(x)\leq V_\infty:=\lim\limits_{|x|\rightarrow\infty}V(x)$ a.e. in $x\in\mathbb{R}^3$ and $3 < r < 5$. For the case of $K=1$ and $g(x,u)=|u|^{r-1} u$ with $2 < r \leq 3$, the existence of a ground state solution to system \eqref{0001500} was proved in \cite{2008ZhaoLG}. In \cite{WXG2017}, when $K=1$, $V$ and $g$ are periodic functions in $x$, the authors obtained a ground state solution for system \eqref{0001500}.
A positive ground state solution of system \eqref{0001500} was obtained in \cite{GG2010}, when $V=1$, $g(x,u)=a(x)|u|^{r-1} u$ with $3 < r <5$, where the nonnegative functions $a$ and $K$ satisfy $a(x)\geq a_\infty:=\lim\limits_{|x|\rightarrow\infty}a(x)$ for $\forall x\in\mathbb{R}^3$, $\lim\limits_{|x|\rightarrow\infty}K(x)=0$.
In \cite{2011Alve},
the authors assumed that $K=1$ and $V$ satisfies \\
$(V_0)$ $ 0< V(x)\leq {V_{T}(x)}$ for $\forall x\in\mathbb{R}^3$, and the strict inequality holds on a set of positive measure, $\lim\limits_{|x|\rightarrow\infty}|V(x)-V_{T}(x)|$=0 and $V_T(x+z)=V_T{(x)}$ for all $x\in \mathbb{R} ^{3}$ and $z\in \mathbb{Z}^3.$\\
They obtained the existence of a positive ground state solution to system \eqref{0001500}.
Concerning the asymptotically periodic system \eqref{0001500},
the same result was obtained in\cite{WXG2017,ZXZ2013} under more weaker assumptions on $V$ and $K$ than that in \cite{2011Alve}.
In \cite{WXG2017}, the authors
supposed that $g$ is an asymptotically periodic function, $V$ and $K$ satisfy \\
$(V)$ $ 0 \leq V(x)\leq {V_{T}(x)}\in {L^{\infty}(\mathbb{R} ^{3})} $ and $ V(x)-{V_T(x)}\in {\mathcal F_{0}}$, where
\begin{eqnarray*}
	\mathcal F_0:={\{k(x):\mathrm{for~any~}\varepsilon>0,\mathrm{meas}\:\{x\in B_1(y):|k(x)|\geq\varepsilon\}\rightarrow0 \:as\: |y|\rightarrow\infty\}}	\end{eqnarray*}
and $V_T$ satisfies \:$V_0:=\inf\limits_{x\in \mathbb{R} ^{3}}V_T>0$ and $V_T(x+z)=V_T{(x)}$ for all $x\in \mathbb{R} ^{3}$ and $z\in \mathbb{Z}^3;$\\
$(K)$ $0< K(x)\leq K_{T}(x)\in L^{\infty}(\mathbb{R}^{3})$, $K(x)-K_{T}(x)\in\mathcal{F}_{0}$, and $K_T$ satisfies $K_{0}:=\inf\limits_{x\in\mathbb{R}^{3}}K_{T}>0$ and $K_T(x+z)=K_T{(x)}$ for all $x\in \mathbb{R} ^{3}$ and $z\in \mathbb{Z}^3$.
\\
The existence of a positive solution for asymptotically periodic system similar to \eqref{0001500} was established in \cite{LiuH2016}, which has critical nonlocal term.

We mention that the following quasilinear Schr\"odinger-Poisson system
\begin{eqnarray}
	\left \{\begin{array}{ll}
		\displaystyle -\Delta_p u+|u|^{p-2}u+\lambda\phi |u|^{m-2}u=|u|^{r-2}u&\mathrm{in} \ \mathbb{R}^3,\\
		\displaystyle -\Delta_q \phi = |u|^m&\mathrm{in}\ \mathbb{R}^3
	\end{array}
	\right. \label{1}
\end{eqnarray}
was first studied in \cite{DSW2023}, where $\lambda>0$, $1<p<3,p<r<p^*$, $q$ and $m$ satisfy
 \begin{eqnarray}\label{10005}
 \displaystyle	\max \left\{1,\frac{3p}{5p-3}\right\}<q<3,\quad\displaystyle \max\left\{1,\frac{(q^*-1)p}{q^*}\right\}<m<\frac{(q^*-1)p^*}{q^*} .	\end{eqnarray}
 Other works regarding the quasilinear \eqref{1}, we refer to \cite{DSW2022cpaa,DSW2022jmaa}.

As far as we know,
there is no results on the existence of ground state solutions to the asymptotically periodic Schr\"odinger-Poisson system \eqref{1001}, this is one of the aims for our paper.
In this article, we
assume that $g\in C(\mathbb{R} ^{3}\times\mathbb R,\mathbb R)$ and satisfies the following hypotheses:
\begin{description}
	\item[$(g_1)$] $\lim\limits_{s\rightarrow{0}}\frac{g(x,s)}{|s|^{p-2}s}=0$ uniformly for $x\in \mathbb{R} ^{3}$;
	\item[$(g_2)$] $\lim\limits_{|s|\rightarrow+\infty}\frac{g(x,s)}{|s|^{p^{*}-2}s}=0$ uniformly for $x\in \mathbb{R} ^{3}$;
	\item[$(g_3)$] $\frac{g(x,s)}{|s|^{\frac{qm}{q-1}-1}}$  is nondecreasing  on $ \mathbb R\backslash\{0\}$;
	\item[$(g_4)$]  there exists $g_T\in C(\mathbb{R} ^{3}\times \mathbb R,\mathbb R)$ such that
	\begin{description}
		\item[{\rm (i) }] $g(x,s)\geq g_T(x,s)$ for all $(x,s)\in \mathbb{R} ^{3}\times \mathbb R^+$,  $g(x,s)\leq g_T(x,s)$ for all $(x,s)\in \mathbb{R} ^{3}\times \mathbb R^-,$
		and $g(x,s)-g_T(x,s)\in\mathcal{F},$
		where\\
		${\mathcal F}:={\{\widetilde{g}(x,s):\mathrm{for~any} \;\varepsilon>0, \mathrm{meas}\:\{x\in B_1(y):|\widetilde{g}(x,s)|\geq\varepsilon\}\rightarrow0 \:\mathrm{as~} |y|\rightarrow\infty}
		\\ ~~~~~~~~\mathrm{~uniformly~for~}|s| \mathrm{~bounded}\};$
		\item[{\rm (ii) }] $g_T(x+z,s)=g_T(x,s)$ for all $(x,s)\in \mathbb{R} ^{3}\times \mathbb R$ and $z\in  \mathbb Z^3;$
		\item[{\rm (iii) }] $\lim\limits_{s\to 0}\frac{g_T(x,s)}{|s|^{\frac{qm}{q-1}-2}s}=0$ uniformly for $x\in \mathbb{R} ^{3};$
		\item[{\rm (iv) }] $\frac{g_T(x,s)}{|s|^{\frac{qm}{q-1}-1}}$  is nondecreasing  on $ \mathbb R\backslash\{0\};$
		\item[{\rm (v) }] $\lim\limits_{s\rightarrow{\infty}}\frac{G_T(x,s)}{|s|^{\frac{qm}{q-1}}}=+\infty$ uniformly for $x\in \mathbb{R} ^{3}$, where $G_T(x,s)=\int^{s}_{0}g_T(x,t)dt.$
	\end{description}
\end{description}

We next consider
the following quasilinear Schr\"odinger-Poisson system with critical exponent
\begin{eqnarray} \label{1002}
	\left\{
	\begin{array}{ll}
		-\Delta_pu+V(x)|u|^{p-2}u+K(x)\phi|u|^{m-2}u=g(x,u)+Q(x)|u|^{p^*-2}u&\quad\text{in}\:\mathbb{R}^3,\\-\Delta_q\phi=K(x)|u|^{m}&\quad\text{in}\:\mathbb{R}^3,
	\end{array}
	\right.
\end{eqnarray}
where $V,K, Q$ and $g$ are asymptotically periodic functions in $x$. The non-existence of solutions, positive and sign-changing solutions, ground states of the semilinear system  \eqref{1002} with $p=q=m=2$ have been studied in \cite{AAP2008, DDM2004, WXG2017, 2014LiuZ, ZXZ2014,ZZ2009, ZT2018}.

Azzollini and   Pomponio \cite{AAP2008}  proved  the existence of a ground state solution for the semilinear system  \eqref{1002} with $p=q=m=2$, $K=Q=1$, $g(x,u)=|u|^{r-1} u$ and $3 < r < 5$,
when $V$ is a positive constant, or $V$ is possibly unbounded below satisfying $V(x)\leq V_\infty:=\lim\limits_{|x|\rightarrow\infty}V(x)$ a.e. in $x\in\mathbb{R}^3$. For the case of $p=q=m=2$, $K=1$ and $g(x,u)=\mu|u|^{r-1} u$ with $\mu>0$, when $3 < r < 5$ or $2 < r < 3$ and $\mu$ is sufficiently large, the existence of a positive ground state solution to the semilinear system \eqref{1002} was established in \cite{2014LiuZ} if $V$ satisfies $V(x)\leq V_\infty:=\liminf\limits_{|x|\rightarrow\infty}V(x)$.
In \cite{WXG2017, ZXZ2014},
the authors assumed that $V, K, Q$  and $g$ are asymptotically periodic functions in $x$, and obtained a ground state solution for the semilinear system  \eqref{1002} with $p=q=m=2$. In \cite{DSW2021}, the authors obtained the ground state solutions and nontrivial solutions to the quasilinear system \eqref{1002} for the case of $q=2$, $m=p$, $V=Q=1,$ $g(x,u)=\lambda|u|^{r-1}u$, where $\frac{3}{2}<p<3$ and $\lambda>0$.

Observe that when $K=0$ and $p=2$. System \eqref{1002} becomes to the semilinear Schr\"odinger equation
\begin{eqnarray} \label{000151}
	-\Delta u + V(x)u= g(x,u)+Q(x)|u|^{2^*-2}u&\quad \text{in~} \mathbb{R}^N.
\end{eqnarray}
In \cite{LLT2016}, the authors studied the existence of a positive
ground state solution for equation \eqref{000151}, when $V$ and $g$ are asymptotically periodic functions in $x$ and $Q\in C(\mathbb{R} ^{N},\mathbb R)$  satisfies \\
$(Q_0)$	there exist a function $Q_T\in C(\mathbb{R}^N,\mathbb{R}),$ 1-periodic in $x_i,1\leq i\leq N$, and a point $x_0\in\mathbb{R}^N$ such that $Q-Q_T\in\mathcal{F}_0$,
$Q(x)\geq Q_T(x)>0$ for all $x\in\mathbb{R}^N$ and
\begin{eqnarray*}Q( x) = |Q| _{\infty }+ O( |x-x_0|^\alpha), \quad\text{as } x\to {x_0},~~\alpha=1 \text{~for~} N=3 \text{~and~} \alpha=2 \text{~for~} N\geq4.
\end{eqnarray*}

Motivated by these works mentioned above, a natural question arises: do ground state solutions for asymptotically periodic Schr\"odinger-Poisson system \eqref{1002} exist or not?
In this article, we consider the  ground state solutions for system \eqref{1002}.	
We assume that  $Q\in C(\mathbb{R}^3,\mathbb{R})$ satisfies\\
$(Q)$ there exist a function $Q_T\in C(\mathbb{R}^3,\mathbb{R}),$ 1-periodic in $x_i,1\leq i\leq 3$, and a point $x_0\in\mathbb{R}^3$ such that $Q-Q_T\in\mathcal{F}_0$,
$Q(x)\geq Q_T(x)>0$ for all $x\in\mathbb{R}^3$ and
\begin{eqnarray*}Q( x) = |Q| _{\infty }+ O( |x-x_0|^\beta) \quad\text{as } x\to {x_0},\quad\text{where } \beta> {0};\end{eqnarray*}
and $g\in C(\mathbb{R}^3,\mathbb{R})$ satisfies ($g_1)$-($g_4)$ and \\
$(g_5)$  $\lim\limits_{s\rightarrow{\infty}}\frac{G(x,s)}{|s|^{\gamma}}=+\infty$ uniformly for $x\in \mathbb{R} ^{3}$, where $G(x,s)=\int^{s}_{0}g(x,t)dt$ and   $\gamma$ satisfies the following hypotheses when	  $1<p\leq\sqrt{3},$
\begin{eqnarray} \label{p*}
	&&p^*>\gamma\geq\left\{
	\begin{array}
		{ll}	\frac{qm}{q-1} ,& \beta\geq\frac{3}{p-1},\\
		\max\left\{\frac{(3-\beta )p}{3-p},\frac{qm}{q-1}\right\},& \beta<\frac{3}{p-1},
	\end{array}\right.
\end{eqnarray} when $\sqrt{3}<p<3,$
\begin{equation}
	\resizebox{0.9\hsize}{!}{$p^*	>\gamma\geq\left\{
		\begin{array}
			{ll}	\max\left\{	{p^*-\frac{mq}{(q-1)(p-1)}},\frac{p(4p-6)}{(3-p)(p-1)} ,\frac{3pq-3p-4q+3}{(3-p)(p-1)}, \frac{qm}{q-1}\right\},&\beta\geq\frac{3}{p-1},\\
			\max\left\{	{p^*-\frac{mq}{(q\!-\!1)(p\!-\!1)}},\frac{p(4p\!-\!6)}{(3\!-\!p)(p-1)} ,\frac{3pq-3p-4q+3}{(3-p)(p-1)}, \frac{(3-\beta )p}{3-p},\frac{qm}{q-1}\right\},&\beta<\frac{3}{p-1}.
		\end{array}\right.$}
\end{equation}

We point out that it is hard to give an explicit expression for the solution of the quasilinear Poisson equation
$-\Delta_q\phi=K(x)|u|^{m}$, which leads to the ineffectiveness of the method used in \cite{WXG2017,WZGcritical2019,ZXZ2013,ZXZ2014} for the proof of Lemma 2.4 in this paper. We have overcome the obstacle by some basic inequalities and calculations of some delicacy.
It is noted that the range of values for $q$ and $m$ in this article is smaller than that in \cite{DSW2023}, since $q$ and $m$ not only need to satisfy \eqref{10005},
but also need to satisfy $\frac{mq}{q-1}<p^*$, which guarantees that $G(x,s)$ is subcritical growth in $s$ at infinity.

Hereafter, we always assume that $1<p<3, q$ and $m$ satisfy \eqref{102}.
The main results are  the  following.
\begin{thm} \label{lemc1.1}
	Assume that $(V),(K)$ and (g$_1)$-(g$_4)$  hold.	Then system \eqref{1001} admits a ground state solution.
\end{thm}
\begin{cor} \label{cor1} Let  $V(x)\equiv V _T(x)$, $K(x)\equiv K_T(x)$  and g$(x, s)$ $\equiv$ g$_T(x, s)$. Assume that $(V),(K),$ (g$_1)$-(g$_4)$  hold. Then system \eqref{1001} admits  a  ground state solution.
\end{cor}
\begin{thm} \label{lemsub1.3}	Assume that $(V),(K),(Q)$ and (g$_1)$-(g$_5)$ hold. Then system \eqref{1002} admits  a  ground state solution.
\end{thm}\begin{cor} \label{cor2} 	Let $V(x)\equiv V_T(x)$, $K(x)\equiv K_T(x)$, $Q(x)\equiv Q_T(x)$ and $g(x, s)\equiv$ g$_T(x, s)$. Assume that $(V),(K),(Q)$, (g$_1)$-(g$_5)$ hold. Then system \eqref{1002} admits a  ground state solution.
\end{cor}
\begin{rem}\label{zhouqifc}	
	{\rm (i) }  In some sense, this paper extend the results in  \cite{2011Alve,WXG2017,WZGcritical2019,ZXZ2013,ZXZ2014} to the case of $p,q\neq2.$ \\
	{\rm (ii) } The proofs of Corollaries 1.2 and 1.4 can be obtained by arguing as in Theorems 1.1 and 1.3, respectively.\\
	{\rm (iii) } The functional sets $\mathcal{F}_0$ in $(V)$ and $\mathcal{F}$ in $(g_4)$ were first introduced in \cite{LLT2016}.
\end{rem}

 The paper is organized as follows. Section 2 is devoted to some notations and essential lemmas. In Section 3 we establish the variational setting of the subcritical system \eqref{1001} and give the proof of Theorem \ref{lemc1.1}. Finally, in Section 4 we deal with the critical system \eqref{1002}, proving Theorem \ref{lemsub1.3}.
	\section{Preliminaries}
	~~~~In this section, we present some  notations and lemmas, which are crucial in proving our results.
	\begin{itemize}
		\item  $W^{1,p}(\mathbb{R}^3)$ is the usual Sobolev space endowed with the norm
		\begin{eqnarray*}\|u\|_{W^{1,p}}=\left(\int_{\mathbb{R}^3}|\nabla u|^p+|u|^p dx\right)^{1/p}.	\end{eqnarray*}
		\item $D^{1,q}(\mathbb{R}^3)$ denotes the closure of  $C_0^\infty(\mathbb{R}^3)$ with respect to the norm $$\|u\|_{D^{1,q}}= \left(\int_{\mathbb{R}^3} |\nabla u|^q dx\right)^{1/q}.$$
		\item $L^{r}(\mathbb{R}^3)$  $(1\leq r<\infty)$ is the Lebesgue space endowed with the norm
		\begin{eqnarray*}|u|_{L^r}=\left(\int_{\mathbb{R}^3} |u|^r dx\right)^{1/r},\end{eqnarray*}
		and	$L^{\infty}(\mathbb{R}^3)$ is the space of measurable functions with the norm
		\begin{eqnarray*}
			|u|_{\infty}=\mathrm{ess}\:\mathrm{sup}\:|u(x)|.
		\end{eqnarray*}
		\item $ E:=\{u\in L^{p^*}(\mathbb{R}^3):|\nabla u|\in L^p(\mathbb{R}^3)$ and $ \displaystyle\int_{\mathbb{R}^3}V(x)|u|^pdx<\infty\}$ is equipped  with the norm
		\begin{eqnarray*}\|u\|=\left(\int_{\mathbb{R}^3}(|\nabla u|^p+V(x)|u|^p)dx\right)^\frac{1}{p}.	\end{eqnarray*}
	\end{itemize}

	We denote by  $E^{*}$ the dual space of the space $E$, which is endowed with the norm $\|\cdot\|_{*}$. $C$ denotes a positive constant that can different from line to line. The notation $o(1)$ represents the quantity that tends to zero.

 Arguing as in
  \cite{LLT2016}, we prove the following lemma.
	\begin{lem} \label{lemc2.2}
	Let $V$ satisfy $(V)$. Then there exist positive constants $c_{1}$ and $c_{2}$ such that for every $u\in E$,
		\begin{eqnarray*}c_{1}\|u\|_{W^{1,p}}^{p}\leq\|u\|^{p}\leq c_{2}\|u\|_{W^{1,p}}^{p}.
		\end{eqnarray*}
		Furthermore, the embedding $E\hookrightarrow  L^l(\mathbb{R}^3)$ is continuous for all $l\in[p,p^*]$. In particular,  $E\hookrightarrow  L_{\text{loc}}^{l}(\mathbb{R}^3)$ is compact if $l\in[p,p^*).$
	\end{lem}
	\textbf{Proof.}  Clearly, the assumption $(V)$ implies  that
	$\|u\|^{p}\leq\max\{1,|V_{T}|_{\infty}\}\|u\|_{W^{1,p}}^{p}$. Then, we  set $c_{2}=\max\{1,|V_{T}|_{\infty}\}$. Let 	\begin{eqnarray*}W(x):=V(x)-V_{T}(x).\end{eqnarray*}
	We deduce from the assumption  $(V)$  that $W(x)\in\mathcal{F}_{0}$. We first show that for any $\varepsilon>0$, there exists $G_{\varepsilon}>0$  such that
	\begin{equation}\label{2.c7}
		\int_{\{x\in\mathbb{R}^3:|W(x)|\geq\varepsilon\}}|u|^{p}dx\leq C\int_{B_{{G_\varepsilon}+1}}|u|^{p}dx+C\varepsilon^{\frac{p}{3}}\|u\|_{W^{1,p}}^{p}\quad\text{for~every}\:u\in E,
	\end{equation}
where  $B_{{G_\varepsilon}+1}:=\{x\in \mathbb{R}^3:|x|<{{G_\varepsilon}+1}\}$, and $C$ is a positive constant that is independent of $\varepsilon.$  Actually,  for any $\varepsilon>0$, there exists $G_\varepsilon>0$ such that for every $|y|\geq G_{\varepsilon}$,
	\begin{equation}\label{aaa}
		\operatorname* { meas} \{ x\in B_{1}( y) : |W(x)| \geq \varepsilon \} < \varepsilon,
	\end{equation}
where $B_{1}(y):=\{x\in \mathbb{R}^3:|x-y|<{1}\}$.
Let  $\mathbb{R}^3$ be covered by balls  $B_1(y_i)$ for $i\in\mathbb{N}$, in this way that every point $x\in\mathbb{R}^3$ lies within at most 4 balls.   For simplicity, we assume that $|y_i|<G_\varepsilon$ for $ i=1,2,3,\ldots,n_\varepsilon$, and $|y_i|\geq G_\varepsilon$ for $i=n_\varepsilon+1,n_\varepsilon+2,n_\varepsilon+3,\ldots,+\infty.$  Applying \eqref{aaa}, the  H\"older and Sobolev inequalities, we obtain that
	\begin{eqnarray*}
		\int_{\{x\in\mathbb{R}^3:|W(x)|\geq \varepsilon\}}|u|^{p}dx
		&\leq&\sum_{i=1}^{+\infty}\int_{\{{x\in }B_{1}(y_{i}):|W(x)|\geq \varepsilon\}}|u|^{p}dx\\
		&=&\sum_{i=1}^{n_{\varepsilon}}\int_{\{{x\in }B_{1}(y_{i}):|W(x)|\geq \varepsilon\}}|u|^{p}dx+\sum_{i=n_{\varepsilon}+1}^{+\infty}\int_{\{{x\in }B_{1}(y_{i}):|W(x)|\geq \varepsilon\}}|u|^{p}dx\\
		&\leq&\!\!\!\!4\int_{\{{x\in}B_{{G_\varepsilon}+1}\!:|W(x)|\geq \varepsilon\}}|u|^{p}dx\!+\!\!\!\!\sum_{i=n_{\varepsilon}+1}^{+\infty}(\text{meas}\{x\!\in\! B_{1}({y}_{i})\!:\!|W(x)|\!\geq\!\varepsilon\})^{\frac{p}{3}}\\
		&&\cdot\left(\int_{\{{x\in}B_{1}({y}_{i}):|W(x)|\geq \varepsilon\}}|u|^{p^{*}}dx\right)^{\frac{3-p}{3}}\\
		&\leq& 4\int_{B_{{G_\varepsilon}+1}}|u|^{p}dx+\left(\sup_{|y|\geq G_{\varepsilon}}\left(\text{meas}\{x\in B_{1}(y):|W(x)|\geq\varepsilon\}\right)\right)^{\frac{p}{3}}\\
		&&\cdot\sum_{i=n_{\varepsilon}+1}^{+\infty}\left(\int_{\{x\in B_{1}(y_i):|W(x)|\geq \varepsilon\}}|u|^{p^{*}}dx\right)^{\frac{3-p}{3}}\\
		&\leq& 4\int_{B_{{G_\varepsilon}+1}}|u|^{p}dx+\varepsilon^{\frac{p}{3}}\sum_{i=1}^{+\infty}\biggl(\int_{B_{1}(y_{i})}|u|^{p^{*}}dx\biggr)^{\frac{3-p}{3}} \\
		&\leq& 4\int_{B_{{G_\varepsilon}+1}}|u|^{p}dx+C\varepsilon^{\frac{p}{3}}\sum_{i=1}^{+\infty}\int_{B_{1}(y_{\mathrm{i}})}(|\nabla u|^{p}+|u|^{p})dx \\
		&\leq& C\int_{B_{{G_\varepsilon}+1}} |u|^{p}dx+C\varepsilon^{\frac{p}{3}}\|u\|_{W^{1,p}}^{p}.
	\end{eqnarray*}
	Thus, the inequality \eqref{2.c7} holds. From $(V)$ we can see that   $|W(x)|\leq 2|V_T|_\infty$. Together with \eqref{2.c7},  by the  H\"{o}lder and Sobolev inequalities, we have
		\begin{eqnarray} \label{bbbbb}
		\left.
		\begin{aligned}
		\int_{\mathbb{R}^{3}}|W(x)||u|^{p}dx&=\int_{\{x\in\mathbb{R}^3:|W(x)|\geq\varepsilon\}}|W(x)||u|^{p}dx+\int_{\{x\in\mathbb{R}^3:|W(x)|<\varepsilon\}}|W(x)||u|^{p}dx \\&\leq C|V_{T}|_{\infty}\int_{B_{G_{\varepsilon}+1}}|u|^{p}dx+C|V_{T}|_{\infty}\varepsilon^{\frac{p}{3}}\|u\|_{W^{1,p}}^{p}+\varepsilon\int_{\mathbb{R}^{3}}|u|^{p}dx \\
			&\leq\! C|V_{T}|_{\infty}\left( \text{meas}B_{{G_\varepsilon}\!+\!1}\right)^{\frac{p}{3}}\left(\int_{\mathbb{R}^{3}}\!\!|u|^{p^{*}}\!\!dx\right)^{\frac{3\!-\!p}{3}}\!\!\!\!\!\!\!+\!C|V_{T}|_{\infty}\varepsilon^{\frac{p}{3}}\|u\|_{W^{1,p}}^p
\!\!	\!+\!\varepsilon\!\!\int_{\mathbb{R}^{3}}\!\!|u|^{p}dx\\
		&\leq C_{\varepsilon}\int_{\mathbb{R}^{3}}|\nabla u|^{p}dx+C|V_{T}|_{\infty}\varepsilon^{\frac{p}{3}}\|u\|^{p}_{W^{1,p}}+\varepsilon\int_{\mathbb{R}^{3}}|u|^{p}dx,
		\end{aligned}
		\right.
	\end{eqnarray}
where $C_\varepsilon$ is a positive constants that is dependent of $\varepsilon$.
From $(V)$, it is clear that $V(x)\geq0$,  $V_T(x)\geq V_0>0$ and $W(x)\leq0.$ Then by \eqref{bbbbb}, we obtain
	\begin{eqnarray*}(C_{\varepsilon}+1)\|u\|^{p}&\geq&(C_{\varepsilon}+1)\int_{\mathbb{R}^{3}}|\nabla u|^{p}dx+\int_{\mathbb{R}^{3}}V(x)|u|^{p}dx\\&=&(C_{\varepsilon}+1)\int_{\mathbb{R}^{3}}|\nabla u|^{p}dx+\int_{\mathbb{R}^{3}}V_{T}(x)|u|^{p}dx-\int_{\mathbb{R}^{3}}|W(x)||u|^{p}dx\\&\geq&(C_{\varepsilon}+1)\int_{\mathbb{R}^{3}}|\nabla u|^{p}dx+V_{0}\int_{\mathbb{R}^{3}}|u|^{p}dx-C_{\varepsilon}\int_{\mathbb{R}^{3}}|\nabla u|^{p}dx\\&&-C|V_{T}|_{\infty}\varepsilon^{\frac{p}{3}}\|u\|_{W^{1,p}}^{p}-\varepsilon\int_{\mathbb{R}^{3}}|u|^{p}dx\\&=&\int_{\mathbb{R}^{3}}|\nabla u|^{p}dx+(V_{0}-\varepsilon)\int_{\mathbb{R}^{3}}|u|^{p}dx-C|V_{T}|_{\infty}\varepsilon^{\frac{p}{3}}\|u\|_{W^{1,p}}^{p}.\end{eqnarray*}
	For $\varepsilon>0$  small enough, there exists  constant $c_1>0$  such that the inequality $c_1\| u\| _{W^{1,p}}^p\leq \| u\|^p$ holds. Moreover, using the Sobolev embedding theorem,  the continuity of the embedding  $E\hookrightarrow L^l(\mathbb{R}^3)$ for all $ l\in[p,p^*]$ and the compactness of the embedding $E\hookrightarrow L_{\text{loc}}^l(\mathbb{R}^3)$ for all $ l\in[p,p^*)$ can be concluded.
	The proof is complete. $\square$

	\begin{lem}\label{0001.1}
		Assume that $(K)$  holds. Then for every $u\in E$, there exists a unique  $\phi_u\in D^{1,q}(\mathbb{R}^3)$ such that
		\begin{eqnarray} \label{01b1}
			\left.
			\begin{array}{ll}
				-\Delta_{q}\phi=K(x)\left|u\right|^{m}.
			\end{array}
			\right.
		\end{eqnarray}
	\end{lem}
	\textbf{Proof.} Fix $u\in E$ and $v\in D^{1,q}(\mathbb{R}^3)$,  thanks to the H\"{o}lder and Sobolev inequalities, we  obtain
	\begin{eqnarray} \label{01.2}
			\left|\int_{\mathbb{R}^{3}}K(x)\lvert u\rvert^{m}vdx\right|\leqslant C|K|_{\infty}\|u\|^{m}\|v\|_{D^{1,q}}.
	\end{eqnarray}
For $u\in E$,  the linear functional $\mathcal{L}$ is defined by
	\begin{eqnarray*}\mathcal{L}(v)=\int_{\mathbb{R}^{3}}K(x)\lvert u\rvert^{m}v dx,\quad v\in D^{1,q}(\mathbb{R}^{3}).	\end{eqnarray*}
Then \eqref{01.2} implies   that the functional $\mathcal{L}(v)$  is well-defined and continuous on $D^{1,q}(\mathbb{R}^3).$ Indeed, when $q=2$,  the result follows directly from the Lax-Milgram  theorem \cite{HB2011}. When $q\neq2,$ by the Minty-Browder theorem \cite{HB2011},   we obtain the desired conclusion using the same arguments as in  \cite[Proposition 2.1]{DSW2023}. The proof is complete.  $\square$

In the same way as in the proof of  Lemma \ref{0001.1}, for any $u\in E$, we can show that the equation
	\begin{eqnarray} \label{posong2}
		\left.
		\begin{array}{ll}
			-\Delta_q \phi = K_T(x) |u|^m \quad \text{~in~}  \mathbb{R}^3
		\end{array}
		\right.
	\end{eqnarray}
has a unique solution which is denoted by $\tilde{\phi}_u$.
	
	Although the precise form of the solution $\phi_u$ is hard to obtain,
 based on \cite[Proposition 2.2]{DSW2023}, we summarize its key properties as follows.
	
	\begin{lem}\label{lem2.0002}
		Let $(K)$  hold. Then for any $u\in E$, the following properties hold
		\begin{enumerate}
			\item[{\rm (i) }] $\displaystyle \int_{\mathbb{R}^{3}}\biggl(\frac{1}{q}|\nabla\phi_{u}|^{q}\!-\!K(x)|u|^{m}\phi_{u}\biggr)dx=\!\!\!\!\min\limits_{{\phi\in D^{1,q}(\mathbb{R}^{3})}}\int_{\mathbb{R}^{3}}\biggl(\frac{1}{q}|\nabla\phi|^{q}\!-\!K(x)|u|^{m}\phi\biggr)dx$ and $\phi_{u}\geqslant0;$
			\item[{\rm (ii) }] for $t>0,\:\phi_{tu}=t^{\frac{m}{q-1}}\phi_{u}$, and for all $y\in \mathbb{R}^{3}$, $\phi_{u(\cdot+y)}=\phi_u{(\cdot+y)};$
			\item[{\rm (iii) }] $\|\phi_{u}\|_{D^{1,q}}\leqslant C\|u\|^{\frac{m}{q-1}}$, where C>0 is independent of $u;$
			\item[{\rm (iv) }] if $u_n\rightharpoonup u$ in $E$, then $\phi_{u_n}\rightharpoonup \phi_u $ in $D^{1,q}({\mathbb{R}^{3}})$ and \begin{eqnarray*}\int_{\mathbb{R}^3}K(x)\phi_{u_n} |u_n|^{m-2}u_n\varphi dx\rightarrow \int_{\mathbb{R}^3}K(x)\phi_u|u|^{m-2}u\varphi dx\quad\text{for all~}\:\varphi\in D^{1,q}(\mathbb{R}^{3}).\end{eqnarray*}
		\end{enumerate}
	\end{lem}

	\begin{lem}\label{wdb zm}
	Assume that $(K)$ holds. If $u_n \rightharpoonup 0$ in $E$, then up to subsequences, we have
	\begin{eqnarray*} \displaystyle \int_{\mathbb{R}^3}\left(K(x)\phi_{u_n} |u_n|^{m-2}u_n\varphi(\cdot-z_n)-K_T(x)\widetilde{\phi}_{u_n} |u_n|^{m-2}u_n\varphi(\cdot-z_n)\right)dx=o{(1)},	\end{eqnarray*}
	where $|z_n|\rightarrow\infty$, $\varphi \in C_0^\infty(\mathbb{R}^3)$, $\phi_{u_n}$ and   $\widetilde{\phi}_{u_n}$ are the  solutions to equations \eqref{01b1} and \eqref{posong2}, respectively.
\end{lem}
\textbf{Proof.}	 We set $H(x):= K(x)-K_T(x)$.  It follows from the assumption $(V)$ that $H(x)\in{\mathcal{F}_0}$. Then
\begin{eqnarray*}
	\begin{array}{ll}
		& \displaystyle \int_{\mathbb{R}^{3}}\left(K(x)\phi_{u_{n}}|u_{n}|^{m-2}u_n\varphi(\cdot-z_{n})-K_{T}(x)\widetilde{\phi}_{u_{n}}|u_{n}|^{m-2}u_n\varphi(\cdot-z_{n})\right)dx\\
		[3mm] =&\displaystyle \int_{\mathbb{R}^{3}}\left((K(x)\!-\!K_T(x))\phi_{u_{n}}|u_{n}|^{m-2}u_n\varphi(\cdot\!-\!z_{n})dx\!+\!\int_{\mathbb{R}^{3}}K_{T}(x){\phi}_{u_{n}}|u_{n}|^{m-2}u_n\varphi(\cdot\!-\!z_{n})\right)dx\\
		[3mm] 	&\displaystyle -\int_{\mathbb{R}^{3}}K_{T}(x)\widetilde{\phi}_{u_{n}}|u_{n}|^{m-2}u_n\varphi(\cdot-z_{n})dx\\
		[3mm] =&\displaystyle \int_{\mathbb{R}^{3}}H(x)\phi_{u_{n}}|u_{n}|^{m-2}u_n\varphi(\cdot-z_{n})dx+\int_{\mathbb{R}^{3}}K_{T}(x)({\phi}_{u_{n}}-\widetilde{\phi}_{u_{n}})|u_{n}|^{m-2}u_n\varphi(\cdot-z_{n})dx\\
		[3mm]  :=&A+B.
	\end{array}
\end{eqnarray*}
To end the proof, it suffices to show that 	\begin{eqnarray}\label{A+B}
	A\rightarrow0,B\rightarrow0\quad \text{as~} n\rightarrow\infty.\end{eqnarray} From the assumption $(K)$, for all $\varepsilon > 0$, there exists  $G_\varepsilon>0$ such that
\begin{eqnarray}\label{Gy}
	\text{meas}\{x \in B_1(y_i): |H(x)| \geq \varepsilon\}<\varepsilon\quad \text{for any }|y|\geq G_\varepsilon.	
\end{eqnarray}
According to  Lemma \ref{lemc2.2}, we know that  $\mathbb{R}^3$ is covered by balls  $B_1(y_i)$ for $i\in\mathbb{N}$, in this way that every point $x\in\mathbb{R}^3$ lies within at most 4 balls. We also assume that $|y_i|<G_\varepsilon$ for $i = 1, 2, ..., n_\varepsilon$, and $|y_i| \geq G_\varepsilon$ for $ i = n_\varepsilon + 1$, $n_\varepsilon + 2$, $n_\varepsilon + 3$, . . . , +$\infty$. Then
\begin{eqnarray}	\left.
	\begin{array}{ll}
		|A|&\leq\displaystyle\int_{\{x\in\mathbb{R}{^3}:|H(x)|\geq\varepsilon\}}\left|H(x)\phi_{u_{n}}|u_{n}|^{m-2}u_n\varphi(\cdot-z_{n})\right|dx\\[3mm]&\displaystyle~~~~+\int_{\{x\in\mathbb{R}{^3}:|H(x)|<\varepsilon\}}\left|H(x)\phi_{u_{n}}|u_{n}|^{m-2}u_n\varphi(\cdot-z_{n})\right|dx\\[3mm]
		&:=A_1+A_2.
	\end{array}
	\right.\label{A=}
\end{eqnarray}
Moreover,  we deduce that
\begin{eqnarray*}
	A_{1} & =&\int_{\{x\in\mathbb{R}{^3}:|H(x)|\geq\varepsilon\}}\left|H(x)\phi_{u_{n}}|u_{n}|^{m-2}u_n\varphi(\cdot-z_{n})\right|dx \\
	&=&\int_{\{x\in\mathbb{R}{^3}:|H(x)|\geq\varepsilon,|x|> G_{\varepsilon}+1\}}\left|H(x)\phi_{u_{n}}|u_{n}|^{m-2}u_n\varphi(\cdot-z_{n})\right|dx \\
	&&+\int_{\{x\in\mathbb{R}{^3}:|H(x)|\geq\varepsilon,|x|\leq G_{\varepsilon}+1\}}\left|H(x)\phi_{u_{n}}|u_{n}|^{m-2}u_n\varphi(\cdot-z_{n})\right|dx \\
	& \leq&\sum_{i=n_{\varepsilon}+1}^{+\infty}\int_{\{x\in B_{1}(y_{i}):|H(x)|\geq\varepsilon,|x|>G_{\varepsilon}+1\}}\left|H(x)\phi_{u_{n}}|u_{n}|^{m-2}u_n\varphi(\cdot-z_{n})\right|dx \\
	&&+\int_{B_{G_{\varepsilon}+1}}\left|H(x)\phi_{u_{n}}|u_{n}|^{m-2}u_n\varphi(\cdot-z_{n})\right|dx \\
	& :=&A_{11}+A_{12}.
\end{eqnarray*}
Combining  \eqref{Gy}, Lemma \ref{lem2.0002}(iii), the H\"older and   Sobolev inequalities, we get
\begin{eqnarray*}
	A_{11} \!\!\!\!& =&\!\!\!\!\!\sum_{i=n_{\varepsilon}+1}^{+\infty}\int_{\{x\in B_{1}(y_{i}):|H(x)|\geq\varepsilon,|x|>G_{\varepsilon}+1\}}\left|H(x)\phi_{u_{n}}|u_{n}|^{m-2}u_n\varphi(\cdot-z_{n})\right|dx  \\
	& \leq &\!\!\!\!\!\sum_{i=n_{\varepsilon}+1}^{+\infty}\!\!\!\! \left(\text{meas}\{x\in B_{1}(y):|H(x)|\geq\varepsilon\}\right)^\frac{{1} }{\lambda}\!\!\cdot\!\left(\int_{\{x\in B_{1}(y_{i}):|H(x)|\geq\varepsilon,|x|>G_{\varepsilon}+1\}}\! \! \! |\phi_{u_n}(x)|^{q^*}dx\right)^\frac{{1} }{q^*}
	\\&&\!\!\!\! \cdot\left(\int_{\{x\in B_{1}(y_{i}):|H(x)|\geq\varepsilon,|x|>G_{\varepsilon}+1\}}\!\! \! \! \! \!|u_{n}|^{p^*}dx\!\!\right)^{\frac{m-1}{p^*}} \!\! \!\!\!\!\cdot\left(\int_{\{x\in B_{1}(y_{i}):|H(x)|\geq\varepsilon,|x|>G_{\varepsilon}+1\}}\!\! \! \!|\varphi(\cdot\!-\!z_{n})|^{m}dx\!\!\right)^{\frac{1}{p^*}}\\
	& \leq&\!\!\!\! C\varepsilon^\frac{{1}}{\lambda}||\phi_{u_n}||_{D^{1,q}}\sum_{i=n_{\varepsilon}+1}^{+\infty}\left(\int_{\{x\in B_{1}(y_{i}):|H(x)|\geq\varepsilon,|x|>G_{\varepsilon}+1\}}|u_{n}|^{p^*}dx\right)^{\frac{m-1}{p^*}} \\
	&&\!\!\!\! \cdot\left(\int_{\{x\in B_{1}(y_{i}):|H(x)|\geq\varepsilon,|x|>G_{\varepsilon}+1\}}|\varphi(\cdot-z_{n})|^{m}dx\right)^{\frac{{1} }{p^*}}\\	& \leq& C\varepsilon^\frac{{1}}{\lambda},\\\text{and}\\
	A_{12} & =&\int_{B_{G_{\varepsilon}+1}}|H(x)\phi_{u_{n}}|u_{n}|^{m-2}u_n\varphi(\cdot-z_{n})|dx \\
	& \leq&   |H|_{\infty}\left(\int_{B_{G_{\varepsilon}+1}}|\phi_{u_{n}}|^{q^*}dx\right)^{\frac{1}{q^*}} \cdot \left(\int_{B_{G_{\varepsilon}+1}}|{u_{n}}|^{\frac{q^*m}{q^*-1}}dx\right)^{\frac{(q^*-1)^2(m-1)}{(q^*)^2m}}\\&& \cdot \left(\int_{B_{G_{\varepsilon}+1}}|\varphi(\cdot-z_{n})|^{\frac{q^*m}{m+q^*-1}}dx\right)^{\frac{(m+q^*-1)(q^*-1)}{(q^*)^2m}}\\	& \leq& C||u_n||^{\frac{m}{q-1}}   \cdot \left(\int_{B_{G_{\varepsilon}+1}}|{u_{n}}|^{\frac{q^*m}{q^*-1}}dx\right)^{\frac{(q^*-1)^2(m-1)}{(q^*)^2m}},
\end{eqnarray*}
where
\begin{eqnarray}\label{a}
	\lambda= \frac{p^*q^*}{p^*q^*-mq^*-p^*}>0.	\end{eqnarray}	
From \eqref{102}, we see that $p<\frac{q^*m}{q^*-1}<p^*$.	Since $u_{n}\rightharpoonup0$ in $E$, it follows from Lemma \ref{lemc2.2} that $u_{n}\rightarrow0$ in $L^{\frac{q^*m}{q^*-1}}_\text{loc}(\mathbb{R}^3)$. Then we conclude that	\begin{eqnarray}\label{A101}
	A_1 \rightarrow0 \quad\text{~as~}	\varepsilon\rightarrow0, n\to\infty. \end{eqnarray}
Using  the H\"older and   Sobolev inequalities, and by Lemma \ref{lem2.0002}(iii), we obtain
\begin{eqnarray*}
	A_2 & =&\int_{\{x\in\mathbb{R}{^3}:|H(x)|<\varepsilon\}}|H(x)\phi_{u_{n}}|u_{n}|^{m-2}u_n\varphi(\cdot-z_{n})|dx \\
	& \leq&\varepsilon\int_{\mathbb{R}^3}|\phi_{u_{n}}|u_{n}|^{m-2}u_n\varphi(\cdot-z_{n})|dx \\
	& \leq&  \varepsilon \left(\int_{\mathbb{R}^3}|\phi_{u_{n}}|^{q^*}dx\right)^{\frac{1}{q^*}}\cdot \left(\int_{\mathbb{R}^3}|{u_{n}}|^{\frac{q^*m}{q^*-1}}dx\right)^{\frac{(q^*-1)^2(m-1)}{(q^*)^2m}}\\&&\cdot \left(\int_{\mathbb{R}^3}|\varphi(\cdot-z_{n})|^{\frac{q^*m}{m+q^*-1}}dx\right)^{\frac{(m+q^*-1)(q^*-1)}{(q^*)^2m}}\\
	& \leq& C \varepsilon||u_n||^{\frac{m}{q-1}}   \cdot||u_n||^{\frac{(q^*-1)(m-1)}{q^*}}
	\\	& \leq&  C \varepsilon.
\end{eqnarray*}
Let $\varepsilon\rightarrow0$, then $A_{2}\rightarrow0$.
Together with \eqref{A=} and \eqref{A101}, we  deduce that 	\begin{eqnarray}\label{A}A\rightarrow0.\end{eqnarray} Furthermore,  thanks to  the H\"older and   Sobolev inequalities, we get
\begin{eqnarray}
	\left.
	\begin{array}{ll}
		|B|&\leq\displaystyle \int_{\mathbb{R}^{3}}|K_{T}(x)({\phi}_{u_{n}}-\widetilde{\phi}_{u_{n}})|u_{n}|^{m-2}u_n\varphi(\cdot-z_{n})|dx\\[3mm]
		&\displaystyle \leq|K_T|_{\infty} \left(\int_{\mathbb{R}^{3}}|{\phi}_{u_{n}}-\widetilde{\phi}_{u_{n}}|^{q^*}dx\right)^{\frac{1}{q^*}} \left(\int_{\mathbb{R}^3}|{u_{n}}|^{\frac{q^*m}{q^*-1}}dx\right)^{\frac{(q^*-1)^2(m-1)}{(q^*)^2m}}\\[3mm]&~~~~\cdot\displaystyle \left(\int_{\mathbb{R}^3}|\varphi(\cdot-z_{n})|^{\frac{q^*m}{m+q^*-1}}dx\right)^{\frac{(m+q^*-1)(q^*-1)}{(q^*)^2m}}
		\\[3mm]
		&\leq  C||{\phi}_{u_{n}}-\widetilde{\phi}_{u_{n}}||_{D^{1,q}}||u_n||^{\frac{(q^*-1)(m-1)}{q^*}}.
	\end{array}
	\right.\label{B1}
\end{eqnarray}
Applying the elementary inequality, there exists $c_q > 0$ such that for any $x,y \in \mathbb{R}^3$,
\begin{eqnarray} \label{jibbds}
	\left\{
	\begin{array}{ll}
		\langle |x|^{q-2} x - |y|^{q-2} y, x-y\rangle_{\mathbb{R}^3 }\ge c_q |x-y|^q&\quad\text{for } 2 \le q < 3,\\(|x| + |y|)^{2-q} \langle |x|^{q-2} x - |y|^{q-2}  y, x-y\rangle_{\mathbb{R}^3 }\ge c_q |x-y|^2&\quad\text{for } 1< q < 2,
	\end{array}
	\right.
\end{eqnarray}
where $\langle\cdot,\cdot\rangle_{\mathbb{R}^3}$ is the standard inner product in $\mathbb{R}^3$.

(i) For $2 \le q < 3$, by \eqref{jibbds}, we have
\begin{eqnarray*}
	C_q||{\phi}_{u_{n}}-\widetilde{\phi}_{u_{n}}||^{q}_{D^{1,q}}\leq \langle-\Delta_q{\phi}_{u_{n}}-(-\Delta_q\widetilde{\phi}_{u_{n}}),{\phi}_{u_{n}}-\widetilde{\phi}_{u_{n}}\rangle.\end{eqnarray*}
Since
\begin{eqnarray*}
	\int_{\mathbb{R}^{3}}-\Delta_q\phi_{u_{n}}(\phi_{u_{n}}-\widetilde{\phi}_{u_{n}})&=&\int_{\mathbb{R}^{3}}K(x)|{u_{n}}|^{m}(\phi_{u_{n}}-\widetilde{\phi}_{u_{n}})dx,\\
	\int_{\mathbb{R}^{3}}-\Delta_q\widetilde\phi_{u_{n}}(\phi_{u_{n}}-\widetilde{\phi}_{u_{n}})&=&\int_{\mathbb{R}^{3}}K_T(x)|{u_{n}}|^{m}(\phi_{u_{n}}-\widetilde{\phi}_{u_{n}})dx,
\end{eqnarray*}
then
\begin{eqnarray}\label{3.22}
	\begin{array}{ll}
		\displaystyle  	C_q||{\phi}_{u_{n}}-\widetilde{\phi}_{u_{n}}||^{q}_{D^{1,q}}&\leq  \displaystyle \langle-\Delta_q{\phi}_{u_{n}},{\phi}_{u_{n}}-\widetilde{\phi}_{u_{n}}\rangle-\langle-\Delta_q\widetilde{\phi}_{u_{n}},{\phi}_{u_{n}}-\widetilde{\phi}_{u_{n}}\rangle\\[3mm]
		&= \displaystyle
		\int_{\mathbb{R}^{3}}(K(x)-K_T(x))|{u_{n}}|^{m}(\phi_{u_{n}}-\widetilde{\phi}_{u_{n}})dx\\	[3mm]&:=D.
	\end{array}
\end{eqnarray}
We claim that  $D\rightarrow0$. Indeed, it is easy to see that
\begin{eqnarray}\label{D0}
	\left.
	\begin{array}{ll}
		|D|&\leq\displaystyle \int_{{\mathbb{R}^{3}}}\left|H(x)|{u_{n}}|^{m}(\phi_{u_{n}}-\widetilde{\phi}_{u_{n}})\right|dx \\[4mm]
		&=\displaystyle\int_{\{x\in\mathbb{R}{^3}:|H(x)|\geq\varepsilon,|x|> G_{\varepsilon}+1\}}\left|H(x)|{u_{n}}|^{m}(\phi_{u_{n}}-\widetilde{\phi}_{u_{n}})\right|dx\\[4mm] &~~~~+\displaystyle\int_{\{x\in\mathbb{R}{^3}:|H(x)|\geq\varepsilon,|x|\leq G_{\varepsilon}+1\}}\left|H(x)|{u_{n}}|^{m}(\phi_{u_{n}}-\widetilde{\phi}_{u_{n}})\right|dx\\[4mm]&~~~~ +\displaystyle\int_{\{x\in{\mathbb{R}^{3}}:|H(x)|<\varepsilon\}}\left|H(x)|{u_{n}}|^{m}(\phi_{u_{n}}-\widetilde{\phi}_{u_{n}})\right|dx \\[4mm]
		&: =D_{1}+D_{2}+D_{3}.
	\end{array}
	\right.
\end{eqnarray}
By \eqref{Gy}, using the H\"older and   Sobolev inequalities, we deduce that
\begin{eqnarray*}
	\left.
	\begin{array}{ll}
		D_{1}&= \!\!\displaystyle\sum\limits_{i=n_{\varepsilon}+1}^{+\infty}\!\!\left(\text{meas}\{x\!\in\!\! B_{1}(y)\!:\!\!|H(x)|\!\geq\!\varepsilon\}\right)^\frac{{1} }{\lambda}
	\!\!	\cdot\!\left(\int_{\{x\in B_{1}(y_{i}):|H(x)|\geq\varepsilon,|x|\!>G_{\varepsilon}+1\}}|{u_n}|^{{p^*}}dx\right)^\frac{{m} }{p^*}
		\\[4mm] & \displaystyle ~~~~    \cdot\left(\int_{\{x\in B_{1}(y_{i}):|H(x)|\geq\varepsilon,|x|>G_{\varepsilon}+1\}}|\phi_{u_{n}}-\widetilde{\phi}_{u_{n}})|^{q^*}dx\right)^\frac{{1} }{q^*}
		\\[4mm]
		&= \displaystyle \varepsilon^\frac{{1} }{\lambda}\sum\limits_{i=n_{\varepsilon}+1}^{+\infty}\left(\int_{\{x\in B_{1}(y_{i}):|H(x)|\geq\varepsilon,|x|>G_{\varepsilon}+1\}}|{u_n}|^{p^*}dx\right)^\frac{{m} }{p^*}
		\\[3mm] & \displaystyle~~~~  \cdot\left(\int_{\{x\in B_{1}(y_{i}):|H(x)|\geq\varepsilon,|x|>G_{\varepsilon}+1\}}|\phi_{u_{n}}-\widetilde{\phi}_{u_{n}})|^{q^*}dx\right)^\frac{{1} }{q^*}
		\\[5mm]
		&\leq \displaystyle C\varepsilon^\frac{{1} }{\lambda}||{u_{n}}||^m ||{\phi}_{u_{n}}-\widetilde{\phi}_{u_{n}}||_{D^{1,q}},
	\end{array}
	\right.\label{D11}
\end{eqnarray*} where ${\lambda}$ is given by \eqref{a}.
Then  we obtain $D_{1}\rightarrow0$ as $\varepsilon\rightarrow0$.  Since $u_n\rightharpoonup0$ in $E$,  we have $u_{n}\rightarrow0$ in $L^{\frac{q^*m}{q^*-1}}_\text{loc}(\mathbb{R}^3)$. The H\"older  inequality yields that
\begin{eqnarray*}	\left.
	\begin{array}{ll}
		D_{2}&=\displaystyle \int_{\{x\in\mathbb{R}{^3}:|H(x)|\geq\varepsilon,|x|\leq G_{\varepsilon}+1\}}\left|H(x)|{u_{n}}|^{m}(\phi_{u_{n}}-\widetilde{\phi}_{u_{n}})\right|dx	\\[3mm]
		&\leq \displaystyle  |H|_{\infty}\left(\int_{B_{{G_\varepsilon}+1}}|{u_n}|^{\frac{q^*m}{q^*-1}}dx\right)^\frac{{q^*-1} }{q^*}
		\cdot\left(\int_{B_{{G_\varepsilon}+1}}|\phi_{u_{n}}-\widetilde{\phi}_{u_{n}})|^{q^*}dx\right)^\frac{{1} }{q^*}
		\\[3mm]&\rightarrow\displaystyle 0,
	\end{array}
	\right.\label{D12}\end{eqnarray*}
as $n\to\infty$.
Moreover, by the H\"older  inequality, we get
\begin{eqnarray*}\label{D2}
	D_3&=&	\int_{\{x\in{\mathbb{R}^{3}}:|H(x)|<\varepsilon\}}\left|H(x)|{u_{n}}|^{m}(\phi_{u_{n}}\!\!\!-\!\widetilde{\phi}_{u_{n}})\right|dx\\
	&\leq&\displaystyle \varepsilon\int_{\mathbb{R}^{3}}\left||{u_{n}}|^{m}(\phi_{u_{n}}-\widetilde{\phi}_{u_{n}})\right|dx\\
	&\leq&\displaystyle  \varepsilon\left(\int_{\mathbb{R}^{3}}|{u_n}|^{\frac{q^*m}{q^*-1}}dx\right)^\frac{{q^*-1} }{q^*}
	\cdot\left(\int_{\mathbb{R}^{3}}|\phi_{u_{n}}-\widetilde{\phi}_{u_{n}})|^{q^*}dx\right)^\frac{{1} }{q^*}\\
	&\leq& C\varepsilon,
\end{eqnarray*}
which implies that 	$D_{2}\rightarrow0$ as $\varepsilon\rightarrow 0.	$
Then, from  \eqref{D0}, we conclude that $D\rightarrow0$. Thus, combining  \eqref{B1} and \eqref{3.22}, we get that for  $2\leq q<3$,
\begin{eqnarray}	\left.
	\begin{array}{ll}
		B\rightarrow0 \quad\text{as~}n\rightarrow \infty.
	\end{array}
	\right.\label{b}\end{eqnarray}

(ii) For $1 < q < 2$, since $\{u_n\}$ is bounded in $E$,  by \eqref{jibbds}, the H\"older  inequality and Lemma \ref{lem2.0002}(iii), we have
\begin{eqnarray*}
	c_q^{\frac{q}{2}}\|\phi_{u_n}-\widetilde{\phi}_{u_n}\|^{q}_{{D}^{1,q}}\!\!\!&\leq&\!\!\!\int_{\mathbb{R}^3}(T(\phi_{u_n},\widetilde{\phi}_{u_n}))^{\frac{q}{2}}(|\nabla \phi_{u_n}|+|\nabla \widetilde{\phi}_{u_n}|)^{\frac{q(2-q)}{2}}dx\\
	&\leq&\!\!\! \langle\!-\!\Delta_q\phi_{u_n}\!-\!(-\!\Delta_ q\widetilde{\phi}_{u_n}),\phi_{u_n}\!-\!\widetilde{\phi}_{u_n}\rangle ^{\frac{q}{2}}\cdot\left(\int_{\mathbb{R}^3}(|\nabla \phi_{u_n}|\!+\!|\nabla \widetilde{\phi}_{u_n}|)^q dx\right)^{\frac{2-q}{2}}	\\
	&\leq&\!\!\! C \langle -\Delta_q\phi_{u_n}-(-\Delta_ q\widetilde{\phi}_{u_n}),\phi_{u_n}-\widetilde{\phi}_{u_n}\rangle  ^{\frac{q}{2}},
\end{eqnarray*}
where  	\begin{eqnarray*}T(\phi_{u_n},\widetilde{\phi}_{u_n})=\langle |\nabla \phi_{u_n}|^{q-2}\nabla \phi_{u_n}-|\nabla \widetilde{\phi}_u|^{q-2}\nabla \widetilde{\phi}_{u_n},\nabla \phi_{u_n}-\nabla \widetilde{\phi}_{u_n}\rangle_{\mathbb{R}^3}.
\end{eqnarray*}
Thus, we obtain
\begin{eqnarray*}	c_q\|\phi_{u_n}-\widetilde{\phi}_{u_n}\|^{2}_{{D}^{1,q}}\leq C \langle -\Delta_q\phi_{u_n}-(-\Delta_ q\widetilde{\phi}_{u_n}),\phi_{u_n}-\widetilde{\phi}_{u_n}\rangle .
\end{eqnarray*}
Using the same argument as in (i), we can derive that
\begin{eqnarray*}	\left.
	\begin{array}{ll}
		\|\phi_{u_n}-\widetilde{\phi}_{u_n}\|_{{D}^{1,q}}\rightarrow0\quad\text{~as~}n\to\infty.
	\end{array}
	\right.\label{c}\end{eqnarray*}
Together with  \eqref{B1}, this yields that  $1<q<2$.
\begin{eqnarray}	\left.
	\begin{array}{ll}
		B\rightarrow0.
	\end{array}
	\right.\label{B3}\end{eqnarray}
Hence,  by  \eqref{A},  \eqref{b} and \eqref{B3}, we know that   \eqref{A+B}  holds.
The proof is complete. $\square$
	\begin{lem}\label{lems2.1}  Assume that  $(g_1),(g_2)$, and  (i),(iii),(iv) of $(g_4)$ hold. Then for all $\delta>0$, there exist $r_\delta>0$ and $C_\delta>0$, we have
		\begin{enumerate}
			\item[{\rm (i) }]  $0\leq g_{T}(x,s)\leq g(x,s)$\  for all $(x,s)\in\mathbb{R}^{3}\times\mathbb{R}^{+};$\\
			$ g(x,s)\leq g_{T}(x,s)\leq 0$ \  for all \ $(x,s)\in\mathbb{R}^{3}\times\mathbb{R}^{-};$
			\item[{\rm (ii) }]  $|g(x,s)|\leq\delta|s|^{p-1}$  for all $(x,s)\in\mathbb{R} ^{3}\times [-r_\delta,r_\delta];$
			\item[{\rm (iii) }]  $|g(x,s)|\leq\delta|s|^{p-1}+C_{\delta}|s|^{p^{*}-1}$ for all $(x,s)\in\mathbb{R}^{3}\times\mathbb{R};$\\
			$0\leq G_{T}(x,s)\leq G(x,s)\leq\frac{\delta}{p}|s|^{p}+\frac{C_{\delta}}{p^{*}}|s|^{p^{*}}$ for all $(x,s)\in\mathbb{R} ^{3}\times\mathbb R;$
			\item[{\rm (iv) }]  $|g(x,s)|\leq C_{\delta}|s|^{p-1}+\delta|s|^{p^{*}-1}$ for all  $(x,s)\in\mathbb{R}^{3}\times\mathbb{R};$
			\item[{\rm (v) }]  for any $\alpha\in(p,p^*)$, there exists  a  constant $C_{\delta,\alpha}>0$   such that\\
			$|g(x,s)|\leq\delta(|s|^{p-1}+|s|^{p^{*}-1})+C_{\delta,\alpha}|s|^{{\alpha}-1}$   for all $(x,s)\in\mathbb{R}^{3}\times\mathbb{R};$
		\\
			$0\leq G_{T}(x,s)\leq G(x,s)\leq\frac{\delta}{p}|s|^{p}+\frac{\delta}{p^{*}}|s|^{p^{*}}+\frac{C_{\delta,\alpha}}{\alpha}|s|^{\alpha}$ for all $(x,s)\in\mathbb{R}^3\times\mathbb{R}.$
		\end{enumerate}
	\end{lem}
	\textbf{Proof.}		The proof  proceeds exactly  as in \cite{LLT2016}, so we omit it here.   $\square$
	
By Lemma \ref{lems2.1}, arguing as in \cite{SL2010}, we conclude the following lemma.
	\begin{lem}\label{lems2.1.1}  Assume that $(g_3)$ and  (i),(iii),(iv) of $(g_4)$ hold. Then
			\begin{eqnarray*}\frac{q-1}{qm}g( x, s) s\geq G( x, s) \geq 0 ~for ~all ~( x, s) \in \mathbb{R} ^{3}\times \mathbb{R}.	\end{eqnarray*}
	\end{lem}

		\begin{lem}\label{lemggGG}
		Assume that   $(g_1)$, $(g_2)$, and (i),(iii),(iv) of $(g_4)$ hold.  Suppose $u_n\rightharpoonup u$ in $E$, then, up to a subsequence, we have
		\begin{eqnarray}
			&&	\int_{\mathbb{R}^{3}}g(x,u_{n})u_{n} dx\to\int_{\mathbb{R}^{3}}g(x,u) udx,\label{Gg1}\\[2mm]
			&&
			\int_{\mathbb{R}^{3}}G(x,u_{n}) dx\to\int_{\mathbb{R}^{3}}G(x,u) dx,\label{Gg2}
			\\[2mm]
			&&
			\int_{\mathbb{R}^{3}}g_T(x,u_{n})u_{n} dx\to\int_{\mathbb{R}^{3}}g_T(x,u) udx,\label{Gg3}
			\\[2mm]
			&&
			\int_{\mathbb{R}^{3}}G_T(x,u_{n}) dx\to\int_{\mathbb{R}^{3}}G_T(x,u) dx.\label{Gg4}
		\end{eqnarray}
	\end{lem}
	\textbf{Proof.}  Since $u_{n}\rightharpoonup u$ in $E$, then $u_{n}\to u$ in $L_{\text{loc}}^{l}(\mathbb{R}^{3})$ for all $l\in[p,p^*)$, and $u_{n}(x)\to u(x)$ a.e. in $\mathbb{R}^{3}$. Define $v_n:=u_n-u$. According to the   Brézis-Lieb lemma \cite{BL1983}, up to a subsequence,  we deduce that
	\begin{eqnarray}\label{0.20}
		\int_{\mathbb{R}^{3}}g(x,v_{n})v_{n}dx +\int_{\mathbb{R}^{3}}g(x,u) udx=	\int_{\mathbb{R}^{3}}g(x,u_{n})u_{n} dx + o(1) \quad \text{as~} n\to\infty.	\end{eqnarray}
 To conclude the proof of \eqref{Gg1}, it only remains
 to show that
	\begin{eqnarray}\label{02.21}
		\int_{\mathbb{R}^{3}}g(x,v_{n})v_{n}dx  \to 0 \quad \text{as~}  n\to\infty.\end{eqnarray}
	Since $\{u_n\}$ is bounded in $E$, then sequence $\{v_n\}$ is also bounded in $E$. Therefore, there exists $\xi>0$ such that
	\begin{eqnarray*}
		||v_n||\leq\xi.
	\end{eqnarray*}
From  Lemma \ref{lems2.1}, for any ${\delta}>0$, there exists positive constant $C_{{\delta},\alpha} $ such that
	\begin{eqnarray*}
	\left|\int_{\mathbb{R}^{3}}g(x,v_{n})v_{n}dx \right|
	&\leq& 	 \delta \int_{\mathbb{R}^{3}}|v_{n}|^{p} d x+C_{{\delta},\alpha}  \int_{\mathbb{R}^{3}} |v_{n}|^{\alpha} d x+ \delta \int_{\mathbb{R}^{3}} |v_{n}|^{p^{*}} d x\\
	&\leq& \delta \left( \xi^{p}+ \xi^{p^{*}}\right)+C_{{\delta},\alpha} \xi^{\alpha},
\end{eqnarray*}
where  $p<\alpha<p^*$ and $C_{{\delta},\alpha} $ depends on ${\delta},\alpha$.   Letting $\delta\to0$ and $n\to\infty$, we obtain that    \eqref{02.21} holds. From \eqref{0.20} and \eqref{02.21}, the relation \eqref{Gg1} is true.  Similarly, \eqref{Gg2}-\eqref{Gg4} follow by the same arguments.   The proof is complete. $\square$

The proof of the next lemma is inspired by \cite{LLT2016}. We include its proof for the sake of completeness.
	\begin{lem} \label{lemc2.4}
	Assume that $(g_1),(g_2)$, and (i),(iii),(iv) of $(g_4)$ hold.  If $\{u_n\}$ is bounded in $E$ and $u_n\to0$ in $L_\mathrm{loc}^l(\mathbb{R}^3)$ for any $l\in[p,p^{*}),$ then,  up to subsequences, we have
	\end{lem}
	\begin{equation}\label{2.81}
	\int_{\mathbb{R}^3}(g(x,u_n)-g_T(x,u_n))u_ndx=o(1),
\end{equation}
	\begin{equation}\label{2.8}
		\int_{\mathbb{R}^3}(G(x,u_n)-G_T(x,u_n))dx=o(1).
	\end{equation}
	\textbf{Proof.}  Let  $\widetilde{g}(x,s):=g(x,s)-g_T(x,s)$, from (i) of $(g_4)$, we see that $\widetilde{g}(x,s)\in{\mathcal{F}}$.  For any $\varepsilon>0$,  there exists   $G_\varepsilon>0$ such that for every $|y|\geq G_\varepsilon$ and $|s|\leq\frac{1}{\varepsilon}$,
	\begin{eqnarray}\label{sssssss}\quad\text{meas}\{x\in B_1(y):|\widetilde{g}(x,s)|\geq\varepsilon\}<\varepsilon.\end{eqnarray}
 Let $\mathbb{R}^3$ be covered by balls  $B_1(y_i)$ for $i\in\mathbb{N}$, in this way that every point $x\in\mathbb{R}^3$ lies within at most 4 balls.  We also assume that $|y_i|<G_\varepsilon$ for $ i=1,2,3,\ldots,n_\varepsilon$, and $|y_i|\geq G_\varepsilon$ for $i=n_\varepsilon+1,n_\varepsilon+2,n_\varepsilon+3,\ldots,+\infty.$                                                                                                               Then, it is easy to deduce that
	\begin{eqnarray*}	\left.
		\begin{array}{ll} 	&~~~~\displaystyle\left|\int_{\mathbb{R}^3}(g(x,u_n)-g_T(x,u_n))u_ndx\right|\\[3mm] &\leq\displaystyle\int_{\mathbb{R}^{3}}|\widetilde{g}(x,u_{n})||u_{n}|dx\\	[3mm]
	&\leq\displaystyle\sum_{i=1}^{n_{\varepsilon}}\int_{B_{1}(y_{i})}|\widetilde{g}(x,u_{n})||u_{n}|dx\!+\sum_{i=n_{\varepsilon}+1}^{+\infty}\int_{B_{1}(y_{i})}|\widetilde{g}(x,u_{n})||u_{n}|dx\\	[3mm]
	&:=\displaystyle J_1+J_2,
		\end{array}
\right.\label{J_1+J_2}
\end{eqnarray*}
where
	\begin{eqnarray*}
			J_1=\sum_{i=1}^{n_{\varepsilon}}\int_{B_{1}(y_{i})}|\widetilde{g}(x,u_{n})|\:|u_{n}|dx,
\end{eqnarray*}
	\begin{eqnarray*}
	\left.
	\begin{array}{ll} J_2
		&	=\displaystyle\sum_{i=n_{\varepsilon}+1}^{+\infty}\int_{\{x\in B_{1}(y_{i}):|\widetilde{g}(x,u_{n})|<\varepsilon\}}|\widetilde{g}(x,u_{n})|\:|u_{n}|dx \\
		&~~~~+\displaystyle\sum_{i=n_{\varepsilon}+1}^{+\infty}\int_{\left\{x\in B_{1}(y_{i}):|u_{n}|\leq\frac{1}{\varepsilon},|\widetilde{g}(x,u_{n})|\geq\varepsilon\right\}}|\widetilde{g}(x,u_{n})|\:|u_{n}|dx \\
		&~~~~\displaystyle+\sum_{i=n_{\varepsilon}+1}^{+\infty}\int_{\left\{x\in B_{1}(y_{i}):|u_{n}|>\frac{1}{\varepsilon},|\widetilde{g}(x,u_{n})|\geq\varepsilon\right\}}|\widetilde{g}(x,u_{n})||u_{n}|dx \\
		&:=\displaystyle J_{21}+J_{22}+J_{23}.	\end{array}
	\right.\label{J_21+J_2}
\end{eqnarray*}
Using Lemma \ref{lems2.1} and the Sobolev inequality, since $u_n\to 0$ in $L_\mathrm{loc}^p(\mathbb{R}^3)$, we derive that
		\begin{eqnarray*}
		J_{1}
		&\leq&\displaystyle4\int_{B_{G_{\varepsilon}+1}} |\widetilde{g}(x,u_{n})| |u_{n}|dx\\
		&\leq&\displaystyle4C_{\delta}\int_{B_{G_{\varepsilon}+1}}|u_{n}|^{p}dx+4\delta\int_{B_{G_{\varepsilon}+1}}|u_{n}|^{p^{*}}dx\\
		&\leq&\displaystyle C\delta+o(1),	
\end{eqnarray*}
which implies that $J_1\to 0$  as $\delta\to 0$ and  $n\to \infty$.
	Applying   Lemma \ref{lems2.1}, \eqref{sssssss}, and by the H\"older and Sobolev inequalities, we have
		\begin{eqnarray*}
			J_{21}&=&\displaystyle\sum_{i=n_{\varepsilon}+1}^{+\infty}\int_{\{x\in B_{1}(y_{i}):|\widetilde{g}(x,u_{n})|<\varepsilon,|u_{n}|\leq r_{\delta}\}}|\widetilde{g}(x,u_{n})| |u_{n}|dx \\
			&&\displaystyle+\sum_{i=n_{\varepsilon}+1}^{+\infty}\int_{\{x\in B_{1}(y_{i}):|\widetilde{g}(x,u_{n})|<\varepsilon,|u_{n}|>r_{\delta}\}}|\widetilde{g}(x,u_{n})| |u_{n}|dx \\
			&\leq&\displaystyle\delta\sum_{i=n_{\varepsilon}+1}^{+\infty}\int_{\{x\in B_{1}(y_{i}):|\widetilde{g}(x,u_{n})|<\varepsilon,|u_{n}|\leq r_{\delta}\}}|u_{n}|^{p}dx\\
			&&\displaystyle+\frac{\varepsilon}{r^{p-1}_{\delta}}\sum_{i=n_{\varepsilon}+1}^{+\infty}\int_{\{x\in B_{1}(y_{i}):|\widetilde{g}(x,u_{n})|<\varepsilon,|u_{n}|>r_{\delta}\}}|u_{n}|^{p}dx \\
			&\leq&\displaystyle 4\delta\int_{\mathbb{R}^{3}}|u_{n}|^{p}dx+\frac{4 \varepsilon}{r^{p-1}_{\delta}}\int_{\mathbb{R}^{3}}|u_{n}|^{p}dx \\
			&\leq& \displaystyle C\delta+\widetilde{C_\delta}\varepsilon,\\
		J_{22}&\leq&\displaystyle\sum_{i=n_{\varepsilon}+1}^{+\infty}\int_{\left\{x\in B_{1}(y_{i}):|u_{n}|\leq\frac{1}{\varepsilon},|\widetilde{g}(x,u_{n})|\geq\varepsilon\right\}}(C_{\delta}|u_{n}|^{p}+\delta|u_{n}|^{p^{*}})dx \\
	&\leq&\displaystyle C_{\delta}\sum_{i=n_{\varepsilon}+1}^{+\infty}\left(\mathrm{meas}\left\{x\in B_{1}\left(y_{i}\right):\left|u_{n}\right|\leq\frac{1}{\varepsilon},\left|\widetilde{g}(x,u_{n})\right|\geq\varepsilon\right\}\right)^{\frac{p}{3}}\\&&\displaystyle
\cdot\left(\int_{B_{1}(y_{i})}|u_{n}|^{p^{*}}dx\right)^{\frac{3-p}{3}}+4\delta\int_{\mathbb{R}^{3}}|u_{n}|^{p^{*}}dx\\&\leq&\displaystyle\varepsilon^{\frac{p}{3}}C_{\delta}\sum_{i=n_{\varepsilon}+1}^{+\infty}\int_{B_{1}(y_{i})}(|\nabla u_{n}|^{p}+|u_{n}|^{p})dx+C\delta\\&\leq&\displaystyle4\varepsilon^{\frac{p}{3}}C_{\delta}\int_{\mathbb{R}^{3}}(|\nabla u_{n}|^{p}+|u_{n}|^{p})dx+C\delta\\&\leq&\displaystyle\varepsilon^{\frac{p}{3}}CC_{\delta}+C\delta,
\end{eqnarray*}
and
		\begin{eqnarray*}
		\left.
		\begin{array}{ll}
			J_{23}&\leq\displaystyle\sum_{i=n_{\varepsilon}+1}^{+\infty}\int_{\left\{x\in B_{1}(y_{i}):|u_{n}|>\frac{1}{\varepsilon}\right\}}(\delta|u_{n}|^{p}+\delta|u_{n}|^{p^{*}}+C_{\delta,\alpha}|u_{n}|^{\alpha})dx	\\
			&\leq\displaystyle4\delta\int_{\mathbb{R}^{3}}(|u_{n}|^{p}+|u_{n}|^{p^{*}})dx+C_{\delta,\alpha}\varepsilon^{p^{*}-\alpha}\sum_{i=n_{\varepsilon}+1}^{+\infty}\int_{\left\{x\in B_{1}(y_{i}):|u_{n}|>\frac{1}{\varepsilon}\right\}}|u_{n}|^{p^{*}}dx\\[3mm]
		&\leq\displaystyle4\delta\int_{\mathbb{R}^{3}}(|u_{n}|^{p}+|u_{n}|^{p^{*}})dx+4C_{\delta,\alpha}\varepsilon^{p^{*}-\alpha}\int_{\mathbb{R}^{3}}|u_{n}|^{p^{*}}dx\\
			&\leq\displaystyle C\delta+CC_{\delta,\alpha}\varepsilon^{p^{*}-\alpha},
		\end{array}
		\right.\label{J_21+J_22}
	\end{eqnarray*}
	where $\widetilde{C_\delta}>0$ is a constant depending on $\delta$. We conclude that $J_2\to 0$ as  $\delta\to 0$ and $\varepsilon\to 0$.
	Then \eqref{2.81} holds.  Moreover, the mean value theorem  yields
	\begin{eqnarray*} \displaystyle G(x,u_{n})-G_{T}(x,u_{n})=\left(g(x,t_{n}u_{n})-g_{T}(x,t_{n}u_{n})\right)u_{n}\end{eqnarray*}
	for some $ t_n\in[0,1].$ Repeating the above argument, we can deduce that \eqref{2.8} holds. The proof is complete. $\square$
	\begin{lem} \label{lemc2.5}
		Assume that  $(V),(g_1),(g_2)$, and (i),(iii),(iv) of $(g_4)$ hold. If $\{u_n\}$ is bounded in $E$, then we obtain
		\begin{eqnarray}\label{2.9}
			\int_{\mathbb{R}^3}(V_T(x)-V(x))|u_n|^{p-2}u_n\varphi(\cdot-z_n)dx=o(1),
		\end{eqnarray}
		\begin{eqnarray}\label{2.10}
			\int_{\mathbb{R}^3}(g(x,u_n)-g_T(x,u_n))\varphi(\cdot-z_n)dx=o(1),
		\end{eqnarray}
	where $|z_n|\to+\infty$ and  $\varphi\in C_0^{\infty}(\mathbb{R}^{3})$.
	\end{lem}
	\textbf{Proof.}
	The proof can be proceeded exactly as in \cite[Lemma 2.5]{LLT2016}. We present the proof for the  sake of completeness. It follows  from the
	 assumption $(V)$ that $W(x)\in{\mathcal{F}_0}$, where $W(x)$ is defined in Lemma \ref{lemc2.2}.
Since $\varphi\in C_0^{\infty}(\mathbb{R}^{3})$, applying the Lebesgue dominated convergence theorem, we deduce that
		\begin{eqnarray} \label{sub 2.12}
		\int_{B_{G_{\varepsilon}+1}}|\varphi(\cdot-z_{n})|^{p}dx=o(1).
	\end{eqnarray}	
	 By  the H\"older  inequality,  \eqref{2.c7}  and \eqref{sub 2.12}, we obtain that
	\begin{equation*}
		\begin{array}{ll}
			& \displaystyle \int_{\mathbb{R}^{3}}\!\left|W(x)|u_{n}|^{p-2}u_n\varphi(\cdot-z_{n})\right|dx \\[3mm]
		= &\!\!\!\displaystyle \int_{\{x\in\mathbb{R}^{3}:|W(x)|\geq\varepsilon\}}\!\!\left|W(x)|u_{n}|^{p-2}u_n\varphi(\cdot\!-\!z_{n})\right|dx\!+\!\!\int_{\{x\in\mathbb{R}^{3}:|W(x)|<\varepsilon\}}\!\!\left|W(x)|u_{n}|^{p-2}u_n\varphi(\cdot\!-\!z_{n})\right|\!dx
			\\[3mm]
			\leq & \displaystyle 2|V_{T}|_{\infty}|u_{n}|^{p-1}_{L^p}\left(\int_{\{x\in\mathbb{R}^{3}:|W(x)|\geq\varepsilon\}}|\varphi(\cdot-z_{n})|^{p}dx\right)^{\frac{1}{p}}+\varepsilon|u_{n}|^{p-1}_{L^p}|\varphi|_{L^p}
			\\[3mm]
			\leq & \displaystyle C\left(C\int_{B_{G_{{\varepsilon}}+1}}|\varphi(\cdot-z_{n})|^{p}dx+C\varepsilon^{\frac{p}{3}}\|\varphi\|_{W^{1,p}}^{p}\right)^{\frac{1}{p}}+C\varepsilon
			\\[5mm]	\leq & \displaystyle
			C\varepsilon^{\frac{1}{3}}+C\varepsilon+o(1)
			.\end{array}
	\end{equation*}
As $\varepsilon\to0$, \eqref{2.9} holds.
Moreover, from  (i) of $(g_4)$,  we have $\widetilde{g}(x,s)\in\mathcal{F}$, where $\widetilde{g}(x,s)$ is defined in Lemma \ref{lemc2.4}.  As the previous lemmas,  $\mathbb{R}^3$ is covered by balls  $B_1(y_i)$. Thus, we obtain that
		\begin{eqnarray*}
				&&\displaystyle\left|\int_{\mathbb{R}^3}(g(x,u_n)-g_T(x,u_n))\varphi(\cdot-z_n)dx\right|\\&\leq&\displaystyle\int_{\mathbb{R}^{3}}|\widetilde{g}(x,u_{n})|\:|\varphi(\cdot-z_{n})|dx\\&\leq&\displaystyle\sum_{i=1}^{n_{\varepsilon}}\int_{B_{1}(y_{i})}|\widetilde{g}(x,u_{n})|\:|\varphi(\cdot\!-\!z_{n})|dx\!+\!\sum_{i=n_{\varepsilon}+1}^{\!+\!\infty}\int_{B_{1}(y_{i})}|\widetilde{g}(x,u_{n})|\:|\varphi(\cdot-z_{n})|dx\\&:=&\displaystyle M_1+M_2.
	\end{eqnarray*}
	Using Lemma \ref{lems2.1} and  the  H\"older  inequality, we have
	\begin{eqnarray*}	\left.
	\begin{array}{ll}
		M_{1}&\leq\displaystyle 4\int_{B_{G_{\varepsilon}+1}}|\widetilde{g}(x,u_{n})| |\varphi(\cdot-z_{n})|dx \\[4mm]
		&\leq\displaystyle 4C_{\delta}\int_{B_{G_{\varepsilon}+1}}|u_{n}|^{p-1} |\varphi(\cdot-z_{n})|dx+4\delta\int_{B_{G_{\varepsilon}+1}}|u_{n}|^{p^{*}-1}|\varphi(\cdot-z_{n})|dx \\
		&\leq\displaystyle  C_{\delta}|u_{n}|^{p-1}_{L^p}\biggl(\int_{B_{G_{\varepsilon}+1}}|\varphi(\cdot-z_{n})|^{p}dx\biggr)^{\frac{1}{p}}+4\delta|u_{n}|_{L^{p^{*}}}^{p^{*}-1}|\varphi|_{L^{p^{*}}} \\
		&\leq\displaystyle  C\delta+o(1),
	\end{array}
	\right.\label{M_1+M_2}
\end{eqnarray*}
which implies that  $M_1\to 0$  as  $\delta\to 0$  and $n\to \infty$. Then, to prove \eqref{2.10}, it suffices to check that
\begin{eqnarray}\label{M20}
	M_{2}\to 0.
\end{eqnarray}
Indeed, it is easy to see that
\begin{eqnarray*}	\left.
	\begin{array}{ll}
		M_{2}
		&=	\displaystyle
		\sum_{i=n_{\varepsilon}+1}^{\!+\!\infty}\int_{\{x\in B_{1}(y_{i}):|\widetilde{g}(x,u_{n})|<\varepsilon\}}|\widetilde{g}(x,u_{n})|\:|\varphi(\cdot-z_{n})|dx\\&	\displaystyle~~~~+
		\sum_{i=n_{\varepsilon}+1}^{+\infty}\int_{\left\{x\in B_{1}(y_{i}):|u_{n}|\leq\frac{1}{\varepsilon},|\widetilde{g}(x,u_{n})|\geq\varepsilon\right\}}|\widetilde{g}(x,u_{n})|\:|\varphi(\cdot-z_{n})|dx\\&~~~~\displaystyle+\sum_{i=n_{\varepsilon}+1}^{+\infty}\int_{\left\{x\in B_{1}(y_{i}):|u_{n}|>\frac{1}{\varepsilon},|\widetilde{g}(x,u_{n})|\geq\varepsilon\right\}}|\widetilde{g}(x,u_{n})|\:|\varphi(\cdot-z_{n})|dx\\&:=\displaystyle M_{21}+M_{22}+M_{23}.
	\end{array}
	\right.\label{M_1+M_2}
\end{eqnarray*}
By Lemma \ref{lems2.1}, \eqref{sssssss}, the  H\"older, Sobolev and Young inequalities, we derive the following estimates
			\begin{eqnarray*}
				M_{21}& \leq&\displaystyle\varepsilon\sum_{i=n_{\varepsilon}+1}^{+\infty}\int_{\{x\in B_{1}(y_{i}):|\widetilde{g}(x,u_{n})|<\varepsilon\}}|\varphi(\cdot-z_{n})|dx \\
				&\leq& \displaystyle4\varepsilon\int_{\mathbb{R}^{3}}|\varphi(\cdot-z_{n})|dx \\
				&\leq&\displaystyle C\varepsilon,\\
			M_{22}&\leq&\displaystyle\sum_{i=n_{\varepsilon}+1}^{+\infty}\int_{\left\{x\in B_{1}(y_{i}):|u_{n}|\leq\frac{1}{\varepsilon},|\widetilde{g}(x,u_{n})|\geq\varepsilon\right\}}(C_{\delta}|u_{n}|^{p-1} |\varphi(\cdot-z_{n})|+\delta|u_{n}|^{p^{*}-1}|\varphi(\cdot-z_{n})|)dx \\
			&\leq&\displaystyle C_\delta\sum_{i=n_\varepsilon+1}^{+\infty}\left(\text{meas}\left\{x\in B_1(y_i):|u_n|\leq\frac{1}{\varepsilon},|\widetilde{g}(x,u_n)|\geq\varepsilon\right\}\right)^{\frac{p}{3}} \\
			&&\displaystyle\cdot\left(\int_{B_{1}(y_{i})}||u_{n}|^{p-1}\varphi(\cdot-z_{n})|^{\frac{3}{3-p}}dx\right)^{\frac{3-p}{3}}+4\delta|u_{n}|_{L^{p^{*}}}^{p^{*}-1}|\varphi|_{L^{p^{*}}} \\
			&\leq&\displaystyle\varepsilon^{\frac{p}{3}}C_{\delta}\sum_{i=n_{\varepsilon}+1}^{+\infty}\biggl[\int_{B_{1}(y_{i})}\biggl(\frac{p-1}{p}{|u_{n}|^{p^*}}+\frac{1}{p}|\varphi(\cdot-z_{n})|^{p^*}\biggr)dx\biggr]^{\frac{3-p}{3}}+C\delta \\
			&\leq&\!\!\displaystyle\varepsilon^{\frac{p}{3}}C_{\delta}C\sum_{i=n_{\varepsilon}+1}^{+\infty}\biggl[\biggl(\frac{p-1}{p}\int_{B_{1}(y_{i})}|u_{n}|^{p^*}dx\biggr)^{\frac{3-p}{3}}\!\!+\!\!\biggl(\frac{1}{p}\int_{B_{1}(y_{i})}|\varphi(\cdot-z_{n})|^{p^*}dx\biggr)^{\frac{3-p}{3}}\biggr]\!\!\!+\!C\delta \\
			&\leq&\displaystyle4\varepsilon^{\frac{p}{3}}C_{\delta}C\left(|u_{n}|_{L^{p^{*}}}^{p}\!+|\varphi|_{L^{p^{*}}}^{p}\right)+C\delta \\
			&\leq&\displaystyle\varepsilon^{\frac{p}{3}}C_{\delta}C+C\delta,\\
\text{and}\\
			M_{23}&\leq&\displaystyle\sum_{i=n_{\varepsilon}+1}^{+\infty}\int_{\left\{x\in B_{1}(y_{i}):|u_{n}|>\frac{1}{\varepsilon}\right\}}(\delta|u_{n}|^{p-1}+\delta|u_{n}|^{p^{*}-1}+C_{\delta,\alpha}|u_{n}|^{\alpha-1})|\varphi(\cdot-z_{n})|dx\\
			&\leq&\displaystyle4\delta\int_{\mathbb{R}^{3}}(|u_{n}|^{p-1} |\varphi(\cdot-z_{n})|+|u_{n}|^{p^{*}-1}|\varphi(\cdot-z_{n})|)dx\\
			&&+\displaystyle C_{\delta,\alpha}\varepsilon^{p^{*}-\alpha}\sum_{i=n_{\varepsilon}+1}^{+\infty}\int_{\left\{x\in B_{1}(y_{i}):|u_{n}|>\frac{1}{\varepsilon}\right\}}|u_{n}|^{p^{*}-1}|\varphi(\cdot-z_{n})|dx\\&
			\leq&\displaystyle4\delta(|u_{n}|^{p-1}_{L^p}|\varphi|_{L^p}+|u_{n}|_{L ^{p^{*}}}^{p^{*}-1}|\varphi|_{L^{p^{*}}})+4C_{\delta,\alpha}\varepsilon^{p^{*}-\alpha}|u_{n}|_{L^{p^{*}}}^{p^{*}-1}|\varphi|_{L^{p^{*}}}
			\\&\leq& \displaystyle C\delta+CC_{\delta,\alpha}\varepsilon^{p^{*}-\alpha}.
	\end{eqnarray*}
	Letting $\varepsilon\to0$ and  $\delta\to0$, we conclude that \eqref{M20} holds.  The proof is complete.  $\square$
		\begin{lem}\label{lemmaV}
	Let $V$ satisfy $(V)$. If $u_n\rightharpoonup 0$ in $E$, then we have
			\begin{eqnarray}\label{12} \int_{\mathbb{R}^{3}}(V(x)-V_T(x))\left|u_n\right|^{p}dx=o(1).
\end{eqnarray}
	\end{lem}
		\textbf{Proof.} Since $u_n\rightharpoonup 0$ in $E$,  we get $u_n\to0$ in $L_{\text{loc}}^l(\mathbb{R}^3)$ for all $l\in[p,p^*).$ It follows from the assumption $(V)$ and  \eqref{2.c7} that
		\begin{eqnarray*} \label{bbbb}
\left|	\int_{\mathbb{R}^{3}}\!(V(x)\!-\!\!V_T(x))\left|u_n\right|^{p}dx\right|\!\!	&\leq&\!\!\!\!\int_{\{x\in\mathbb{R}^3:|W(x)|\geq\varepsilon\}}\!\!|W(x)||u_n|^{p}dx\!+\!\!\int_{\{x\in\mathbb{R}^3:|W(x)|<\varepsilon\}}\!\!|W(x)||u_n|^{p}dx \\
			&\leq&\!\!\!\! C|V_{T}|_{\infty}\!\int_{B_{G_{\varepsilon}+\!1}}|u_n|^{p}dx\!+\!C|V_{T}|_{\infty}\varepsilon^{\frac{p}{3}}\|u_n\|_{W^{1,p}}^{p}\!+\!\varepsilon\!\int_{\mathbb{R}^{3}}\!|u_n|^{p}dx,
	\end{eqnarray*}
where $W(x)$ is defined in Lemma \ref{lemc2.2}. Letting
$\varepsilon\to 0$ and $n\to\infty$, then \eqref{12} holds. The proof is complete. $\square$

	\section{The subcritical case}
~~~~	In this section, we establish the variational framework for studying the solutions of system \eqref{1001} and present the proof of Theorem \ref{lemc1.1}.

	The energy functional for system \eqref{1001} is  defined by
	\begin{eqnarray*}\mathcal{J}(u)=\frac{1}{p}\int_{\mathbb{R}^{3}}(|\nabla u|^p+V(x)|u|^p)dx+\frac{q-1}{qm}\int_{\mathbb{R}^{3}}K(x)\phi_u|u|^mdx-\int_{\mathbb{R}^3}G(x,u)dx.	\end{eqnarray*}
	
		\begin{lem}\label{Lemma3.1}
		Assume that $(V),(K),(g_1)$ and $(g_2)$  hold. Then $\mathcal{J} \in C^{1}(E,\mathbb{R} )$ and for every $u,v\in E,$ we have
		\begin{eqnarray*}\langle\mathcal{J} ^{\prime}(u),v\rangle\!\!=\!\!\int_{\mathbb{R}^{3}}\left(|\nabla u|^{p-2}\nabla u\nabla v\!+\!V(x)|u|^{p-2}uv\right)dx\!\!+\!\!\int_{\mathbb{R}^{3}}K(x)\phi_{u}|u|^{m-2}uvdx\!\!-\!\int_{\mathbb{R}^{3}}g(x,u)vdx.\quad\:\end{eqnarray*}
	\end{lem}
	\textbf{Proof.}  The proof can be obtained working as in \cite[Proposition 2.3]{DSW2023} and we omit the details here.   $\square$
	
	According to \cite[Proposition 2.4]{DSW2023}, $(u,\phi)\in W^{1,p}(\mathbb{R}^{3})\times D^{1,q}(\mathbb{R}^{3})$ is a solution of system \eqref{1001} if and only if $u$ is a critical point of $\mathcal{J} $ and $\phi=\phi_{u}$. Moreover,
	 with the same arguments of \cite[Lemma  3.1]{DSW2023}, we can prove the following lemma.
	 \begin{lem}\label{sjx}
	 	Assume that $(V),(K),(g_1)$ and $(g_2)$  hold. If $\{u_n\}$ is a bounded sequence in $E$ satisfying $\mathcal{J}^\prime(u_n)\to0,$ then there exists $u \in E$ such that, up to a subsequence, $\nabla u_n(x) \to \nabla u(x)$ a.e. in $\mathbb{R}^3$.
	 \end{lem}

	In what follows, we denote by $	\mathcal{N}$ the Nehari manifold associated to the functional $\mathcal{J}$, that
	is,
	\begin{equation*}
		\mathcal{N} := \left\{u \in E \setminus \{0\} : \langle\mathcal{J} ^{\prime}(u),u\rangle = 0\right\} \end{equation*}
and  set \begin{eqnarray}\label{ccc}m := \inf\limits_{u \in \mathcal{N}} \mathcal{J}(u).\end{eqnarray}

	\begin{lem} \label{lemc3.1}
		Assume  that   $(V),(K)$, and $(g_1)$-$(g_4)$ hold. For every $u\in E\backslash\{0\},$ there exists a unique $t_u>0$ such that $t_uu\in\mathcal{N}$ and $\mathcal{J}(t_uu)=\max\limits_{t>0}\mathcal{J}(tu)$.
	\end{lem}
	\textbf{Proof.}
	For fixed  $u\in E\backslash\{0\},$    we define a function  $\ell:\mathbb{R}^+ \to \mathbb{R}$ as $\ell(t)=\mathcal{J}(tu)$.	Applying Lemma \ref{lems2.1}(iii), the Fatou lemma  and   (v) of $(g_4)$, we obtain
	\begin{eqnarray*}
		\operatorname*{lim\inf}_{t\to+\infty}\int_{\mathcal{M}}\frac{G(x,tu)}{|tu|^{\frac{qm}{q-1}}}|u|^{\frac{qm}{q-1}}dx\geq\operatorname*{lim\inf}_{t\to+\infty}\int_{\mathcal{M}}\frac{G_T(x,tu)}{|tu|^{\frac{qm}{q-1}}}|u|^{\frac{qm}{q-1}}dx=+\infty,
	\end{eqnarray*}
	where  $\mathcal{M}:=\{x\in\mathbb{R}^3:u(x)\neq0\}.$ Then, since $p<\frac{qm}{q-1}$, by Lemma \ref{lem2.0002}, we deduce that
		\begin{eqnarray*}\ell(t)
&\leq	& t^{\frac{qm}{q-1}}\left(\frac{t^{p-\frac{qm}{q-1}}}{{p}}\|u\|^{p}+\frac{q-1}{qm}\int _{\mathbb{R}^3} K(x)\phi_{u}|u|^{m}dx-\int_{\mathcal{M}}\:\frac{G(x,tu)}{|tu|^{\frac{qm}{q\!-\!1}}}|u|^{\frac{qm}{q\!-\!1}}dx\right)\\
		&\to&-\infty.\end{eqnarray*}
Using  Lemma \ref{lems2.1}(iii) and   the Sobolev inequality, we get
	\begin{eqnarray}\label{123}
		\int_{\mathbb{R}^{3}}G(x,tu)dx\!\leq\!\frac{\delta}{p}t^{p}\int_{\mathbb{R}^{3}}|u|^{p}dx\!+\!\frac{C_{\delta}}{p^{*}}t^{p^{*}}\int_{\mathbb{R}^{3}}|u|^{p^{*}}dx\!\leq\! C\delta t^{p}\|u\|^{p}\!+\!C_{\delta}Ct^{p^{*}}\|u\|^{p^{*}}.
	\end{eqnarray}
	Then, by  the assumption $(K)$, Lemma \ref{lem2.0002} and \eqref{123},  we derive that
	\begin{eqnarray*}\ell(t)&=&\frac{t^{p}}{p}||u||^{p}+\frac{q-1}{qm}\int_{\mathbb{R}^3} K(x)\phi_{tu}|tu|^{m}dx-\int_{\mathbb{R}^3} G(x,tu)dx\\
		&\geq&\left(\frac{1}{p}-C\delta\right)t^{p}\|u\|^{p}-C_{\delta}Ct^{p^{*}}\|u\|^{p^{*}}.\end{eqnarray*}
Choose $\delta>0$ such that $\frac{1}{p}-C\delta>0$. It follows that $\ell(t_0)>0$ for some $t_0>0$ small enough. Consequently, there exists a $t_u>0$ such that $\mathcal{J}(t_uu)=\max\limits_{t>0}\mathcal{J}(tu)$ and $t_uu\in\mathcal{N}.$ If  $\overline{t}_u\neq t_u$  is positive constant that satisfies   $\overline{t}_uu\in\mathcal{N}$, then
	\begin{eqnarray}\label{5}
		{{\overline{t}_u}^{p-\frac{qm}{q-1}}}\|u\|^p+\int_{\mathbb{R}^3}K(x)\phi_u|u|^mdx=\int_\mathcal{M}\dfrac{g(x,\overline{t}_uu)|u|^{\frac{qm}{q-1}-1}u}{{|\overline{t}_uu|}^{\frac{qm}{q-1}-1} }dx.\end{eqnarray}
	Replacing $\overline{t}_u$ by $t_u$, \eqref{5} still holds. Thus we have
	\begin{eqnarray*}	\left(t_{u}^{p-\frac{qm}{q-1}}-\overline{t}_{u}^{~ p-\frac{qm}{q-1}}\right)\|u\|^{p}=\int_\mathcal{M}\left(\frac{g(x,t_{u}u)}{|{t_{u}u}|^{\frac{qm}{q-1}-1}}-\frac{g(x,\overline{t}_{u}u)}{|\overline{t}_{u}u|^{\frac{qm}{q-1}-1}}\right)|u|^{\frac{qm}{q-1}-1}udx.
	\end{eqnarray*}
	From  $(g_3)$  and $\overline{t}_u\neq t_u$, we reach a contradiction. Thus $t_u$ is unique, and the proof is complete. $\square$
	
	The following remark can be proved working as in  \cite{PHR1992,MW1996}.
	\begin{rem} \label{remark3.2}We have
		\begin{eqnarray*}\inf_{u\in\mathcal{N}}\mathcal{J}(u)=\inf_{u\in E\backslash\{0\}}\max_{t>0}\mathcal{J}(tu)=\inf_{\eta\in\Gamma}\max_{t\in[0,1]}\mathcal{J}(\eta(t))>0,\end{eqnarray*}
		where
		\begin{eqnarray*}\Gamma:=\{\eta\in C([0,1],E):\eta(0)=0,\mathcal{J}(\eta(1))<0\}.\end{eqnarray*}
	\end{rem}
	\begin{lem} \label{lemc3.3}
		Assume that  $(V),(K)$ and $(g_1)$-$(g_4)$ hold.  Then there exists a bounded sequence $\{u_n\}\subset E$ satisfying
			\begin{eqnarray}\label{c0}
			\mathcal{J}(u_n) \to m,\quad \|\mathcal{J}^{\prime}(u_n)\|_{{*}}\to0\quad\text{as~}n\rightarrow \infty.\end{eqnarray}
	\end{lem}
	\textbf{Proof.} According to Lemma \ref{lemc3.1},   $\mathcal{J}$ has the mountain pass geometry. Then the  mountain pass theorem \cite{AAPH1973}   yields the existence of a sequence $\{u_n\} \subset E$ satisfying \eqref{c0}.
	Furthermore, by Lemma \ref{lems2.1.1}, we get
	\begin{eqnarray*}
		m&=&\mathcal{J}(u_n)-\frac{q-1}{qm}\langle \mathcal{J}^{\prime}(u_{n}),u_{n}\rangle+o(1) \\
		& =&\left(\frac{1}{p}-\frac{q-1}{qm}\right)\|u_{n}\|^{p}+\int_{\mathbb{R}^{3}}\left(\frac{q-1}{qm}g(x,u_{n}) u_{n}-G(x,u_{n})\right)dx +o(1)\\
		& \geq&\left(\frac{1}{p}-\frac{q-1}{qm}\right)\|u_{n}\|^{p}+o(1),
	\end{eqnarray*}
which implies that the sequence $\{u_{n}\}$ is bounded in $E$ 	because $\frac{1}{p}-\frac{q-1}{qm}>0$.  The proof is complete.   $\square$
	
	\begin{lem} \label{lemc3.4}
		Assume that  $(V),(K)$ and $(g_1)$-$(g_4)$ hold. Let $u\in\mathcal{N} $ and  $\mathcal{J}(u)=m$, then $(u, {\phi}_u)$  is  a solution of system \eqref{1001}.
	\end{lem}
		\textbf{Proof.} The argument to prove this lemma is essentially the same  as that of  \cite{LLT2016},  thus we
		omit the details. $\square$
	
We now present the periodic system that corresponds to the asymptotically periodic system given in \eqref{1001},
	\begin{eqnarray} \label{Vp1001}
		\left\{
		\begin{array}{ll}
			-\Delta_pu+V_T(x)|u|^{p-2}u+K_T(x)\phi|u|^{m-2}u=g_T(x,u)&\quad\text{in}\:\mathbb{R}^3,\\-\Delta_q\phi=K_T(x)|u|^{m}&\quad\text{in}\:\mathbb{R}^3.
		\end{array}
		\right.
	\end{eqnarray}
Let us define the energy functional of  system \eqref{Vp1001} on $E,$	\begin{eqnarray*}\mathcal{J}_T(u)&=&\frac{1}{p}\int_{\mathbb{R}^{3}}(|\nabla u|^p+V_T(x)|u|^p)dx+\frac{q-1}{qm}\int_{\mathbb{R}^{3}} K_T(x)\phi_{u}|u|^{m}dx-\int _{\mathbb{R}^{3}}G_T(x,u)dx,
\end{eqnarray*}and
	\begin{eqnarray*}m_T:=\inf\limits_{u\in \mathcal{N}_{T}} \mathcal{J}_T(u),\end{eqnarray*}
	where
\begin{eqnarray*}\mathcal{N}_{T}=\{u\in E\setminus\{0\}:\langle  \mathcal{J}_{T}'(u),u\rangle=0\}.	\end{eqnarray*}

	\begin{rem}\label{rem 3.5}
 From  Lemma $ \ref{lemc3.1}$,	for each  $u\in E\backslash\{0\}$, there exists $t_u>0$ such that $t_uu\in\mathcal{N}$, which implies that  $\mathcal{J}(t_uu)\geq m.$ Then,  by   the assumptions  $(V),(K)$ and Lemma $\ref{lems2.1}(iii),$ we deduce that	\begin{eqnarray*}m\leq \mathcal{J}(t_uu)\leq \mathcal{J}_T(t_uu)\leq\max\limits_{t>0}\mathcal{J}_T(tu). \end{eqnarray*}
 It follows that
$m\leq \inf\limits_{u\in E\backslash\{0\}} \max\limits_{t>0}\mathcal{J}_T(tu). $
		Under the assumptions on  $g_T$, the  remark $\ref{remark3.2}$ still holds with $\mathcal{J}$ and $\mathcal{N}$ replaced by $\mathcal{J}_T$ and $\mathcal{N}_T$, respectively. Hence,
		 \begin{eqnarray*}m\leq \inf\limits_{u\in E\backslash\{0\}} \max\limits_{t>0}\mathcal{J}_T(tu)= m_T. \end{eqnarray*}
		
	\end{rem}
	\textbf{Proof of Theorem 1.1.} Lemma \ref{lemc3.3} yields the existence of a bounded sequence $\{u_n\}\subset E$ satisfying \eqref{c0}. It follows that there exists  $u\in E$ such that, up to a subsequence, \begin{eqnarray*} \label{E}
		\left.
		\begin{array}{ll}
			\begin{cases}u_n\rightharpoonup u&\text{in}~E,\\u_n\to u&\text{in}~L_{\text{loc}}^l(\mathbb{R}^3)\text{~for all~} l\in[p,p^*),\\u_n(x)\to u(x)&\text{a.e.}\:\text{in~} \mathbb{R}^3.\end{cases}
		\end{array}
		\right.
	\end{eqnarray*}   Then, applying Lemma \ref{lem2.0002}(iv), Lemma \ref{sjx} and  \cite[Proposition 5.4.7]{MW2013}, we deduce that
	$\langle  \mathcal{J}^\prime(u),\varphi\rangle=0$ for every $\varphi\in{C}_0^\infty(\mathbb{R}^3)$. Thus $(u,\phi_u)$ is a solution of system \eqref{1001}. We will prove the existence of a ground state solution for system \eqref{1001}.
	
	Assume $u=0$ and set
	\begin{eqnarray*} \sup_{z\in\mathbb{R}^3}\int_{B_\delta(z)}|u_n|^pdx\to\mu.\end{eqnarray*}
We claim that $ \mu>0$. Indeed, if $ \mu=0$,   from \cite[Lemma 3.2]{DSW2023}, we have $u_n\to0$ in $L^{{l}}(\mathbb{R}^{3})$ for  all $l\in(p,p^*).$ Using $\frac{q-1}{qm}-\frac{1}{p}<0$, the assumption $(K)$,   Lemma \ref{lem2.0002}(i) and  Lemma \ref{lems2.1}(v), we obtain
	\begin{eqnarray*}m&=&\mathcal{J}(u_n)-\frac{1}{p}\langle \mathcal{J}^{\prime}(u_{n}),u_{n}\rangle+o(1)\\&=&\left(\frac{q\!-\!1}{qm}\!-\!\frac{1}{p}\right)\int_{\mathbb{R}^{3}}K(x)\phi_{u_{n}}|u_{n}|^m dx\!-\!\int_{\mathbb{R}^{3}}\left(G(x,u_{n})\!-\!\frac{1}{p}g(x,u_{n}) u_{n}\right)dx\!+\!o(1)\\&\leq&\frac{1}{p}\int_{\mathbb{R}^{3}}g(x,u_{n}) u_{n}dx+o(1)\\&\leq&\frac{1}{p} \int_{\mathbb{R}^{3}}\left(\delta\left(|u_{n}|^{p}\!+\!|u_{n}|^{p^*}\right)\!+\!C_{\delta,\alpha} |u_{n}|^{\alpha}\right)dx+o(1)
		\\&\leq& o(1).
	\end{eqnarray*}
This is in contradiction with $m>0$.  Thus $ \mu>0,$ as we claimed. Consequently,  up to a subsequence,  there exists  a sequence $\{z_n\}\subset\mathbb{Z}^3$ satisfying
	\begin{eqnarray}\label{1a1}
	\int_{B_\delta(z_n)}|u_n|^pdx\geq\frac{\mu}{2}>0.\end{eqnarray}
Let $\overline{u}_n(x):=u_{n}(x                                                                                                                                                                                                                                                                                                                                                                                                                                                         +z_{n}).$ Then, there exists  $\overline{u}\in E$ such that up to a subsequence,
\begin{eqnarray} \label{E1}
	\left.
	\begin{array}{ll}
		\begin{cases}\overline{u}_n\rightharpoonup \overline{u}&\text{in}~E,\\\overline{u}_n\to \overline{u}&\text{in}~L_{\text{loc}}^l(\mathbb{R}^3)\text{~for all~} l\in[p,p^*),\\\overline{u}_n(x)\to \overline{u}(x)&\text{a.e.}\:\text{in~} \mathbb{R}^3.\end{cases}
	\end{array}
	\right.
\end{eqnarray}
It follows from \eqref{1a1} and \eqref{E1} that $\int_{B_\delta}|\overline{u}|^pdx\geq\frac{\mu}{2}>0$, which means that
 $\overline{u}\neq0$. We will show that  $\{z_{n}\}$ is unbounded. Actually, if  $\{z_{n}\}$ is bounded,   then there exists a constant  $G^{\prime}>0$ such that
	\begin{eqnarray*}\int_{B_{G^\prime}}|u_n|^pdx\geq\int_{B_\delta(z_n)}|u_n|^pdx\geq\frac{\mu}{2}.\end{eqnarray*}
This is a contradiction, since  $u_n\to 0$ in $L_{\mathrm{loc}}^l(\mathbb{R}^3)$  for all  $l\in[p,p^*)$. Hence, $\{z_n\}$ is unbounded. Then, up to a subsequence, $|z_n|\to\infty.$
	For every $\varphi\in C_0^{\infty}(\mathbb{R}^{3})$, using  \eqref{E1}, Lemmas \ref{wdb zm} and  \ref{sjx}, by \cite[Proposition 5.4.7]{MW2013},  Lemmas  \ref{lemc2.5} and \ref{lems2.1} , we conclude that
		\begin{eqnarray*}
\langle\mathcal{J}_{T}^{\prime}(\overline{u}),\varphi\rangle
&=&\int_{\mathbb{R}^{3}}\left(|\nabla \overline{u}|^{p-2}\nabla \overline{u}\nabla\varphi+V_T(x)|\overline{u}|^{p-2} \overline{u}\varphi\right) dx+\int_{\mathbb{R}^{3}}K_T(x)\widetilde{\phi}_{\overline{u}}|\overline{u}|^{m-2}\overline{u}\varphi dx
\\&&-\int_{\mathbb{R}^{3}}g_T(x,\overline{u})\varphi dx\\
&=&\!\!\!\int_{\mathbb{R}^{3}}\!\!\left(|\nabla \overline{u}_n|^{p-2}\nabla \overline{u}_n\nabla\varphi\!+\!V_T(x)|\overline{u}_n|^{p-2} \overline{u}_n\varphi\right) dx\!+\!\int_{\mathbb{R}^{3}}K_T(x)\widetilde{\phi}_{\overline{u}_n}|\overline{u}_n|^{m-2}\overline{u}_n\varphi dx\\&&-\int_{\mathbb{R}^{3}}g_T(x,\overline{u}_n)\varphi dx+o(1)\\
	&=&\int_{\mathbb{R}^{3}}\left(|\nabla u_{n}|^{p-2}\nabla u_{n}\nabla\varphi(\cdot-z_{n})+V_T(x)|u_{n}|^{p-2} u_{n}\varphi(\cdot-z_{n})\right)dx\\&&+\int_{\mathbb{R}^{3}}K_T(x)\widetilde{\phi}_{u_n}|u_{n}|^{m-2}u_{n}\varphi(\cdot-z_{n})dx-\int_{\mathbb{R}^{3}}g_T(x,u_{n})\varphi(\cdot-z_{n})dx+o(1) \\
		&=&\int_{\mathbb{R}^{3}}\left(|\nabla u_{n}|^{p-2}\nabla u_{n}\nabla\varphi(\cdot-z_{n})+V(x)|u_{n}|^{p-2} u_{n}\varphi(\cdot-z_{n})\right)dx\\&&+\int_{\mathbb{R}^{3}} K(x)\phi_{u_n}|u_{n}|^{m-2}u_{n}\varphi(\cdot-z_{n})dx-\int_{\mathbb{R}^{3}}g(x,u_{n})\varphi(\cdot-z_{n})dx+o(1)\\
		&=&\langle\mathcal{J}^{\prime}(u_{n}),\varphi(\cdot-z_{n})\rangle+o(1)\\	&=&0.
\end{eqnarray*}
	Thus, $(\overline{u}, \widetilde{\phi}_{\overline{u}})$  is the  solution of system \eqref{Vp1001}. Applying Lemmas \ref{lemggGG}, \ref{lemc2.4} and \ref{lemmaV}, the Fatou lemma,  we deduce that	
			\begin{eqnarray*}
			\left.
			\begin{array}{ll}
		&~~~~m_T\\\!&\displaystyle\leq\mathcal{J}_T(\overline{u})-\frac{q-1}{qm}\langle\mathcal{J}_T^{\prime}(\overline{u}),\overline{u}\rangle\\[3mm]	&\displaystyle=\left(\frac{1}{p}\!-\!\!\frac{q\!-\!1}{qm}\right)\int_{\mathbb{R}^{3}}\left(|\nabla \overline{u}|^{p}\!\!+\!V_T(x)|\overline{u}|^{p} \right) dx\!+\!\int_{\mathbb{R}^{3}}\left(\frac{q\!-\!1}{qm}g_T(x,\overline{u}) \overline{u}\!-\!G_T(x,\overline{u})\right)dx\\[4mm]	&\displaystyle\leq\left(\frac{1}{p}\!-\!\frac{q\!\!-\!1}{qm}\right)\int_{\mathbb{R}^{3}}\!\!\!\left(|\nabla \overline{u}_n|^{p}\!\!+\!\!V_T(x)|\overline{u}_n|^{p} \right) dx\!+\!\!\!\int_{\mathbb{R}^{3}}\!\left(\frac{q\!-\!1}{qm}g_T(x,\overline{u}_n) \overline{u}_n\!\!-\!\!G_T(x,\overline{u}_n)\right)dx\!+o(1)\\[4mm]	&\displaystyle=\!\left(\frac{1}{p}\!-\!\frac{q\!-\!1}{qm}\right)\!\int_{\mathbb{R}^{3}}\!\!\left(|\nabla u_n|^{p}\!\!+\!\!V_T(x)|u_n|^{p} \right) dx\!+\!\int_{\mathbb{R}^{3}}\!\!\left(\!\frac{q\!-\!1}{qm}g_T(x,u_{n}) u_{n}\!\!-\!G_T(x,u_{n})\right) \!dx\!+\!o(1)	\\[4mm]&\displaystyle=\!\left(\frac{1}{p}\!-\!\frac{q\!-\!1}{qm}\right)\!\int_{\mathbb{R}^{3}}\!\left(|\nabla u_n|^{p}\!+\!\!V(x)|u_n|^{p} \right) dx\!+\!\!\int_{\mathbb{R}^{3}}\!\left(\!\frac{q\!-\!1}{qm}g(x,u_{n}) u_{n}\!\!-\!G(x,u_{n})\right) \!dx\!+\!o(1)	\\[3mm]&\displaystyle=\mathcal{J}(u_{n})-\frac{q-1}{qm}\langle\mathcal{J}^{\prime}(u_{n}),u_{n}\rangle+o(1)\\[1mm]&\displaystyle=m.
		\end{array}
	\right.
	\end{eqnarray*}
From  Remark \ref{rem 3.5}, we see  that $\mathcal{J}_{T}(\overline{u})=m_{T}=m>0.$  It follows that $\overline{u}\neq0$. According to   Lemma \ref{lemc3.1}, there exists a unique $t_{\overline{u}}>0$ such that $t_{\overline{u}}\overline{u}\in\mathcal{N}.$ Then, it is easy to check   that
	\begin{eqnarray*}\label{J00}m\leq \mathcal{J}(t_{\overline{u}}\overline{u})\leq \mathcal{J}_T(t_{\overline{u}}\overline{u})\leq \mathcal{J}_T(\overline{u})=m.\end{eqnarray*}
	Thus $\mathcal{J}(t_{\overline{u}}\overline{u})=m$. By Lemma \ref{lemc3.4}, we conclude that ($t_{\overline{u}}\overline{u},\phi_{t_{\overline{u}}\overline{u}})$ is the ground solution of system \eqref{1001}.
	
Assume $u\neq0.$ We know from \eqref{ccc} that	\begin{eqnarray}\label{J00}	
		\mathcal{J}(u)\geq m.\end{eqnarray} Using the assumption $(V)$, the Fatou  lemma and Lemma \ref{lemggGG}, we obtain that
	\begin{eqnarray}
		\left.
		\begin{array}{ll}
			\mathcal{J}(u)&=\displaystyle\mathcal{J}(u)-\frac{q-1}{qm}\langle \mathcal{J}^{\prime}(u),u\rangle\\[3mm]
			&=\displaystyle\left(\frac{1}{p}-\frac{q-1}{qm}\right)||u||^{p}+\int_{\mathbb{R}^{3}}\left(\frac{q-1}{qm}g(x,u) u-G(x,u)\right)dx\\[4mm]
			&\leq\displaystyle\left(\frac{1}{p}-\frac{q-1}{qm}\right)\|u_{n}\|^{p}+\int_{\mathbb{R}^{3}}\left(\frac{q-1}{qm}g(x,u_{n}) u_{n}-G(x,u_{n})\right)dx+o(1)\\[3mm]
			&=\displaystyle\mathcal{J}(u_{n})-\frac{q-1}{qm}\langle \mathcal{J}^{\prime}(u_{n}),u_{n}\rangle+o(1)\\[2mm]
			&=m.	\end{array}
		\right.\label{D120}\end{eqnarray}
	Then, we deduce from \eqref{J00} and  \eqref{D120} that $\mathcal{J}(u)=m$.
Therefore, $(u,{\phi}_u)$  is the ground state solution of system \eqref{1001}. The proof is complete.  $\square$
	\section{The critical case}
~~~~In this section, we give the proof of Theorem \ref{lemsub1.3}.
	Let us define the functional $\mathcal{\widetilde J}$ on $E$ by
	\begin{eqnarray*}\mathcal{\widetilde  J}(u)&=&\frac{1}{p}\int_{\mathbb{R}^{3}}(|\nabla u|^p+V(x)|u|^p)dx+\frac{q-1}{qm}\int_{\mathbb{R}^{3}}K(x)\phi_u|u|^mdx-\int_{\mathbb{R}^3}G(x,u)dx\\&&-\frac{1}{p^*}\int_{\mathbb{R}^{3}}Q(x)|u|^{p^{*}}dx.\end{eqnarray*}
 By Lemma \ref{Lemma3.1}, $\mathcal{\widetilde J}(u)\in C^1(E,\mathbb{R})$ and for any $u, v\in E$,
	\begin{eqnarray*}
		\langle\mathcal{\widetilde J^{\prime}}(u),v\rangle
		&=&\int_{\mathbb{R}^{3}}\left(|\nabla u|^{p-2}\nabla u\nabla v+V(x)|u|^{p-2} uv\right)dx+\int_{\mathbb{R}^{3}}K(x)\phi_{u}|u|^{m-2}uvdx\\&&-\int_{\mathbb{R}^{3}}g(x,u)vdx-\int_{\mathbb{R}^{3}}Q(x)|u|^{p^{*}-2}uvdx.
	\end{eqnarray*}
Clearly, if $u$ is a  critical point of $\mathcal{\widetilde J}(u)$, then
$(u, {\phi}_u)$ is a solution  of system \eqref{1002}.	
	\begin{lem}\label{sjx2}
		Assume  that  $(V),(K),(Q),(g_1)$ and $(g_2)$  hold. Let  $\{u_n\}$ be a bounded sequence in $E$   satisfying $\widetilde {\mathcal{J}}^\prime(u_n)\to0.$ Then, there exists $ u \in E$ such that, up to a subsequence, $\nabla u_{n}( x) \to \nabla u( x)$ a.e. in $ \mathbb{R} ^{3}.$
	\end{lem}
	\textbf{Proof.} The technique is the same as in \cite[Lemma 2.2]{DSW2022cpaa} for this proof and  we omit the details.  $\square$

	\begin{lem}\label{lem2.4}
	Let $Q$ satisfy  $(Q)$ and $\{u_n\} \subset E$ be bounded. Then,
		\begin{eqnarray}\label{c2.011}
			\int_{\mathbb{R}^3}(Q(x)-Q_T(x))|u_n|^{p^*-2}u_n\varphi(\cdot-z_n)dx=o(1),
		\end{eqnarray}
where  $|z_n|\to\infty$ and  $\varphi\in C_0^{\infty}(\mathbb{R}^{3})$.
	\end{lem}
	\textbf{Proof.} The proof proceeds  as in  \cite[Lemma 2.5]{LLT2016}.
We define $\overline{Q}(x):=Q(x)-Q_{T}(x)$. By the assumption $(Q),$ we  see that $\overline{Q}(x)\in\mathcal{F}_0$. Using \eqref{2.c7} and \eqref{sub 2.12}, for every $\varepsilon>0$, there exists $G_{\varepsilon}>0$ such that
		\begin{eqnarray}
	\left.
	\begin{array}{ll}
	\displaystyle\int_{\{x\in\mathbb{R}^3:|\overline{Q}(x)|\geq\varepsilon\}}|\varphi(\cdot-z_{n})|^{p}dx&	\displaystyle\leq C\int_{B_{G_\varepsilon+1}}|\varphi(\cdot-z_{n})|^{p}dx+C\varepsilon^{\frac{p}{3}}\|\varphi\|_{W^{1,p}}^{p}\\&\leq C\varepsilon^{\frac{p}{3}}+o(1).
\end{array}
\right.\label{abb}\end{eqnarray}
	Then, thanks to   the  H\"older  inequality and \eqref{abb}, we deduce that
	\begin{eqnarray*} 	&&\left|\int_{\mathbb{R}^3}(Q(x)-Q_T(x))|u_n|^{p^*-2}u_n\varphi(\cdot-z_n)dx\right|\\&\leq&\int_{\{x\in\mathbb{R}^3:|\overline{Q}(x)|\geq\varepsilon\}}|Q(x)-Q_{T}(x)|\:|u_{n}|^{p^{*}-1}|\varphi(\cdot-z_{n})|dx\\&&+\int_{\{x\in\mathbb{R}^3:|\overline{Q}(x)|<\varepsilon\}}\left|Q(x)-Q_{T}(x)\right|\:|u_{n}|^{p^{*}-1}\left|\varphi(\cdot-z_{n})\right|dx\\
		&\leq&2|Q|_{\infty}|u_{n}|_{L^{p^{*}}}^{p^{*}-1}|\varphi|_{\infty}^{\frac{p^{*}-p}{p^{*}}}\left(\int_{\{x\in\mathbb{R}^3:|\overline{Q}(x)|\geq\varepsilon\}}|\varphi(\cdot-z_{n})|^{p}dx\right)^{\frac{1}{p^{*}}}+\varepsilon|u_{n}|_{L^{p^{*}}}^{p^{*}-1}|\varphi|_{L^{p^{*}}}\\
		&\leq &C\varepsilon^{\frac{3-p}{9}}+C\varepsilon+o(1).
	\end{eqnarray*}
Letting $\varepsilon\to0$,   we conclude that \eqref{c2.011} holds.  The proof is complete. $\square$
	
	The Nehari manifold $\widetilde{\mathcal{N}}$ corresponding to $\widetilde{\mathcal{J}}$ is defined by
	\begin{eqnarray*}\widetilde{\mathcal{N}}=\{u\in E\setminus\{0\}:\langle  \widetilde{\mathcal{J}}'(u),u\rangle=0\},	\end{eqnarray*}
	and set
	\begin{eqnarray*}\label{c00}\widetilde{m}:=\inf\limits_{u\in \widetilde{\mathcal{N}}} \widetilde{\mathcal{J}}(u).\end{eqnarray*}
	
	\begin{lem}\label{lem0003.1}
		Assume that  $(V),(K),(Q)$ and $(g_1)$-$(g_4)$  hold.  Then for every $u\in E\backslash\{0\}$, there exists a  unique $\widetilde{t}_u>0$ such that $\widetilde{t}_uu\in\widetilde{\mathcal{N}}$ and  $\widetilde{\mathcal{J}}(\widetilde{t}_uu)=\max\limits_{t>0}\widetilde{\mathcal{J}}(tu)$.
	\end{lem}
	{\bf Proof.}
We fix $u\in E$ with $u \neq0$, and we consider the function $\widetilde{\ell}(t):=\mathcal{\widetilde J}(tu)$, where $t>0$.
Arguing as we have done in the proof of Lemma \ref{lemc3.1}, it is easy to see that  	\begin{eqnarray*}\widetilde{\ell}(t)\to -\infty\quad \text{~as~}t\to+\infty.\end{eqnarray*}
By  Lemma \ref{lems2.1}(iii) and the Sobolev inequality, we deduce that
	\begin{eqnarray}\label{9000}\int_{\mathbb{R}^{3}}G(x,tu)dx\leq\frac{\delta}{p}t^{p}\int_{\mathbb{R}^{3}}|u|^{p}dx+\frac{C_{\delta}}{p^{*}}t^{p^{*}}\int_{\mathbb{R}^{3}}|u|^{p^{*}}dx\leq C\delta t^{p}\|u\|^{p}+CC_{\delta}t^{p^{*}}\|u\|^{p^{*}}.	\end{eqnarray}
	Then, from  the assumptions $(K), (Q),$   Lemma \ref{lem2.0002} and  \eqref{9000}, we have
	\begin{eqnarray*}\widetilde{\ell}(t)&=&\frac{1}{p}||tu||^p+\frac{q-1}{qm}\int_{\mathbb{R}^{3}} K(x)\phi_{tu}|tu|^{m}dx-\int_{\mathbb{R}^{3}} G(x,tu)dx-\frac{1}{p^*}\int_{\mathbb{R}^{3}} Q(x)|tu|^{p^*}dx\\
		&\geq&\frac{t^{p}}{p}||u||^{p}-\int_{\mathbb{R}^{3}} G(x,tu)dx-\frac{t^{p^*}}{p^*}|Q|_{\infty}\int_{\mathbb{R}^{3}} |u|^{p^*}dx\\
		&\geq&\left(\frac{1}{p}-C\delta\right)t^{p}\|u\|^{p}-CC_{\delta}t^{p^{*}}\|u\|^{p^{*}}.	\end{eqnarray*}
	 Since $\delta$ is arbitrary and $p<p^*$, this shows that $\widetilde{\ell}(t)>0$ for small  $t>0$.
Consequently,  there exists  a $\widetilde{t}_u>0$ such that $\mathcal{\widetilde J}(\widetilde{t}_uu)=\max\limits_{t>0}\mathcal{\widetilde J}(tu)$  and $\widetilde{t}_uu\in\widetilde{\mathcal{N}}$. We claim that $\widetilde{t}_uu$ is unique. Actually, let   $ \overline{t}_u \neq\widetilde{t}_u$ be positive constant satisfying $\overline{t}_uu\in\widetilde{\mathcal{N}}$. It follows that
	\begin{eqnarray}\label{123a}
	\overline{t}_u^{p-\frac{qm}{q-1}}\|u\|^p\!+\!\int_{\mathbb{R}^3}K(x)\phi_u|u|^mdx\!=\!\!\int_{\mathcal{M}}\dfrac{g(x,\overline{t}_uu)|u|^{\frac{qm}{q-1}-1}u}{|\overline{t}_uu|^{\frac{qm}{q-1}\!-\!1} }dx\!+\!\overline{t}_u^{p^*-\frac{qm}{q-1}}\!\int_{\mathbb{R}^{3}} Q(x)|u|^{p^*}dx,\end{eqnarray}
	where ${\mathcal{M}}$ is given by  Lemma \ref{lemc3.1}.	Since \eqref{123a} still hold when $\overline{t}_u$ is replaced by  $\widetilde{t}_u$. Then we obtain
	\begin{eqnarray*}\left(\widetilde{t}_{u}^{p-\frac{qm}{q-1}}-\overline{t}_{u}^{~ p-\frac{qm}{q-1}}\right)\|u\|^{p}&=&
		\int_{\mathcal{M}}\left(\frac{g(x,\widetilde{t}_{u}u)}{|\widetilde{t}_{u}u|^{\frac{qm}{q-1}-1}}-\frac{g(x,\overline{t}_{u}u)}{|\overline{t}_{u}u|^{\frac{qm}{q-1}-1}}\right)|u|^{\frac{qm}{q-1}-1}u dx\\&&+\left(\widetilde{t}_{u}^{p^*-\frac{qm}{q-1}}-\overline{t}_{u}^{~ p^*-\frac{qm}{q-1}}\right)	\int_{\mathbb{R}^{3}} Q(x)|u|^{p^*}dx.\end{eqnarray*}
From $(g_3)$ and $p<\frac{qm}{q-1}<p^*$,  this yields a contradiction with $\overline{t}_u\neq \widetilde{t}_u.$ Thus, the claim follows and the proof is complete. $\square$

With essentially the same arguments as in \cite{PHR1992,MW1996}, the following remark can be proved.
	\begin{rem} \label{sremark3.2}
		We have
		\begin{eqnarray*}\inf_{u\in\widetilde{\mathcal{N}}}\widetilde{\mathcal{J}}(u)=\inf_{u\in E\backslash\{0\}}\max_{t>0}\widetilde{\mathcal{J}}(tu)=\inf_{\widetilde{\eta}\in\widetilde{\Gamma}}\max_{t\in[0,1]}\widetilde{\mathcal{J}}(\widetilde{\eta}(t))>0,\end{eqnarray*}
		where
		\begin{eqnarray*}\widetilde{\Gamma}:=\{\widetilde{\eta}\in C([0,1],E)\colon\widetilde{\eta}(0)=0,\widetilde{\mathcal{J}}(\widetilde{\eta}(1))<0\}.\end{eqnarray*}
	\end{rem}
	\begin{lem}\label{lemsub 3.3} Assume that $(V),(K),(Q)$ and $(g_1)$-$(g_4)$ hold. Then there exists a bounded sequence $\{u_n\}\subset E$ satisfying
			\begin{eqnarray}\label{c01}
				\mathcal{\widetilde J}( u_n) \to  \widetilde{m} ,\quad \| \mathcal{\widetilde J^{\prime }}( u_n) \| _* \to 0 \quad \text{as~}  n\to \infty.\end{eqnarray}
	\end{lem}
	\textbf{Proof.} By  Lemma \ref{lem0003.1} and the  mountain pass theorem \cite{AAPH1973}, we obtain  the existence of  a sequence $\{u_n\} \subset E$ satisfying \eqref{c01}.
Moreover, applying Lemma \ref{lems2.1.1}(i) and  the assumption $(Q)$, we deduce that
	\begin{eqnarray*} \widetilde{m}&\geq &\mathcal{\widetilde J}(u_{n})-\frac{q-1}{qm}\langle \mathcal{\widetilde J}^{\prime}(u_{n}),\:u_{n}\rangle\\
		&=&\left(\frac{1}{p}-\frac{q-1}{qm}\right)\|u_{n}\|^{p}+\int_{\mathbb{R}^{3}}\left(\frac{q-1}{qm}g(x,u_{n}) u_{n}-G(x,u_{n})\right)dx\\&&+\left(\frac{q-1}{qm}-\frac{1}{p^*}\right)\int_{\mathbb{R}^{3}} Q(x)|u|^{p^*}dx \\
		&\geq&\:\left(\frac{1}{p}-\frac{q-1}{qm}\right)\|u_{n}\|^{p}.
	\end{eqnarray*}
It follows  that the sequence $\{u_n\}$ is bounded in $E$,  since $\frac{1}{p}\!-
 \frac{q-1}{qm}>0$. The proof is complete.    $\square$
	\begin{lem} \label{lems4.4}
		Assume that  $(V),(K),(Q)$ and $(g_1)$-$(g_4)$ hold. If $u\in\widetilde{\mathcal{N} }$ and  $\widetilde{\mathcal{J}}(u)=\widetilde{m}$, then $(u, {\phi}_u)$ is  a solution of system \eqref{1002}.
	\end{lem}
{\bf Proof.} 	The proof follows the same argument as in \cite[Lemma 3.4]{LLT2016}, we omit it here.  $\square$

 According to\cite{GT1976},  $S_p$ is the best Sobolev constant for    $D^{1,p}(\mathbb{R}^3)\hookrightarrow L^{p^*}(\mathbb{R}^3)$,  that is,
\begin{eqnarray}\label{Sp}S_p=\displaystyle\inf\limits_{u\in D^{1,p}(\mathbb{R}^3)\backslash\{0\}}\frac{\int_{\mathbb{R}^3}|\nabla u|^pdx}{{\left(\int_{\mathbb{R}^3}|u|^{p^*}dx\right)}^{\frac{p}{p^*}}}.\end{eqnarray}
The achieving function of $S_p$ is denoted by
\begin{eqnarray}\label{2.02}
	U_{\varepsilon,p}(x)=\frac{\varepsilon^{\frac{3-p}{p(p-1)}}\left(\frac{3(3-p)^{p-1}}{(p-1)^{p-1}}\right)^{\frac{3-p}{p^{2}}}}{\left(\varepsilon^{\frac{p}{p-1}}+|x|^{\frac{p}{p-1}}\right)^{\frac{3-p}{p}}},\quad\varepsilon>0.
\end{eqnarray}
Let $\rho>0$ and  $\psi\in C_0^\infty(\mathbb{R}^3)$ satisfy
	\begin{eqnarray*}\label{2.04}
	\left\{\begin{array}{ll}\psi(x)=1&\quad\mathrm{for}\:|x|\leq \rho,\\0\leq\psi(x)\leq1&\quad\mathrm{for}\:\rho<|x|<2\rho,\\\psi(x)=0&\quad\mathrm{for}\:|x|\geq2\rho.\end{array}\right.
\end{eqnarray*}
  We set
\begin{eqnarray}\label{2.03}
	u_{\varepsilon,p}(x)=\psi(x)U_{\varepsilon,p}(x).
\end{eqnarray}
From  \cite{WS2023},
with a simple computation one finds that
\begin{eqnarray}\label{2.05}
	\int_{\mathbb{R}^3}|\nabla u_{\varepsilon,p}|^pdx=S_p^{\frac{3}{p}}+O\left(\varepsilon^{\frac{3-p}{p-1}}\right),
\end{eqnarray}
\begin{eqnarray}\label{2.06}
	\int_{\mathbb{R}^3}| u_{\varepsilon,p}|^{p^*}dx=S_p^{\frac{3}{p}}+O\left(\varepsilon^{\frac{3}{p-1}}\right),
\end{eqnarray}
\begin{eqnarray}\label{2.07}\left.\int_{\mathbb{R}^3}|u_{\varepsilon,p}|^pdx=\left\{
	\begin{array}
		{ll}	O\left(\varepsilon^{\frac{3-p}{p-1}}\right)+b_1\varepsilon^p, & 1<p<\sqrt{3}, \\[2mm]
		O\left(\varepsilon^p\right)+b_2\varepsilon^p|\ln\varepsilon|, & p=\sqrt{3}, \\[2mm]
	O\left(\varepsilon^{\frac{3-p}{p-1}}\right),& \sqrt{3}<p<3,
	\end{array}\right.\right.\quad\end{eqnarray}
\begin{eqnarray}\label{2.08}
	\left.\int_{\mathbb{R}^3}|u_{\varepsilon,p}|^\gamma dx=\left\{
	\begin{array}
		{ll}O\left(\varepsilon^{\frac{\gamma(3-p)}{p(p-1)}}\right),& \gamma<\frac{3(p-1)}{3-p}, \\[2mm]
		O\left(\varepsilon^{\frac{3}{p}}\right)+b_3\varepsilon^{\frac{3}{p}}|\ln\varepsilon| ,& \gamma=\frac{3(p-1)}{3-p}, \\[2mm]
		O\left(\varepsilon^{\frac{\gamma(3-p)}{p(p-1)}}\right)+b_4\varepsilon^{3-\frac{\gamma(3-p)}{p}},& \gamma>\frac{3(p-1)}{3-p},
	\end{array}\right.\right.\quad\end{eqnarray}
\begin{eqnarray}\label{2.07x}\left.\int_{\mathbb{R}^3}|x|^\beta|u_{\varepsilon,p}|^{p^*}dx=\left\{
	\begin{array}
		{ll}	O\left(\varepsilon^{\frac{3}{p-1}}\right)+b_5\varepsilon^\beta,& \beta<\frac{3}{p-1}, \\[2mm]
		O\left(\varepsilon^\frac{3+\beta}{p}\right)+b_6\varepsilon^\frac{3+\beta}{p}|\ln\varepsilon|, &\beta=\frac{3}{p-1}, \\[2mm]
	O\left(\varepsilon^{\frac{3}{p-1}}\right), & \beta>\frac{3}{p-1},
	\end{array}\right.\right.\quad\end{eqnarray}
where $\varepsilon>0$ is small and  $b_i (i=1,2,3,4)$ are positive constants independent of
$\varepsilon.$ For the Poisson term, using \eqref{2.08}, elementary computations yield that
\begin{eqnarray}\label{2.09}
	\int_{\mathbb{R}^3}K(x)\phi_{u_{\varepsilon,p}}|u_{\varepsilon,p}|^mdx\leq\left\{\begin{array}{ll}O\big(\varepsilon^{\frac{mq(3-p)}{p(p-1)(q-1)}}\big)+C\varepsilon^{\frac{4pq-3qm+pqm-3p}{p(q-1)}},&\quad1\!<\!p\!<\!\frac{3qm+4q-3}{qm+4q-3},\\[2mm]O\big(\varepsilon^{\frac{4q-3}{p(q-1)}}\big)+C\varepsilon^{\frac{4q-3}{p(q-1)}}\big|\ln\varepsilon\big|^{\frac{4q-3}{3(q-1)}},&\quad p=\frac{3qm+4q-3}{qm+4q-3},\\[2mm]O\big(\varepsilon^{\frac{mq(3-p)}{p(q-1)(p-1)}}\big),&\quad\frac{3qm+4q-3}{qm+4q-3}\!<\!p\!<\!3,\end{array}\right.\end{eqnarray}
where $\frac{3}{2} < \frac{3qm+4q-3}{qm+4q-3} < 3$ for $p>1$. Actually, by  the assumption $(K)$,  Lemma \ref{lem2.0002}, the  H\"{o}lder  and Sobolev inequalities,   we deduce that
\begin{eqnarray}\label{2.010}
	\int_{\mathbb{R}^3}K(x)\phi_{u_{\varepsilon,p}}|u_{\varepsilon,p}|^mdx\leq C \left(\int_{\mathbb{R}^3}|u_{\varepsilon,p}|^{\frac{mq^*}{q^*{-1}}}dx\right)^{\frac{(q^*-1)q}{q^*(q-1)}}.
\end{eqnarray}
Using \eqref{2.08} and \eqref{2.010}, we obtain the estimate \eqref{2.09}.
\begin{lem}\label{lemsub3.6} Assume that $(V),(K),$ $(Q)$ and $(g_1)$-$(g_5)$ hold. Then the following inequality holds
	\begin{eqnarray}\label{S}
		\inf_{u\in E\setminus\{0\}}\max_{t>0}\mathcal{\widetilde J}(tu)<\frac{1}{3}S^{\frac{3}{p}}|Q|_{\infty}^{-\frac{3-p}{p}}.\end{eqnarray}
\end{lem}
	\textbf{Proof.} To prove \eqref{S}, it is enough to check that
for the function $u_{\varepsilon,p}$ defined in  \eqref{2.03},
\begin{eqnarray*}\label{111}
	\max_{t>0}\mathcal{\widetilde J}(t u_{\varepsilon,p})<\frac{1}{3}S_{p}^{\frac{3}{p}}|Q|_{\infty}^{-\frac{3-p}{p}}\quad\mathrm{as}\:\varepsilon>0\mathrm{~small.}\end{eqnarray*}
Applying \eqref{2.05} and \eqref{2.06},  we deduce that
\begin{eqnarray*}
	& &  \max_{t>0}\left(\frac{t^{p}}{p}\int_{\mathbb{R}^{3}}|\nabla u_{\varepsilon,p}|^{p}dx-\frac{t^{p^{*}}}{p^{*}}\int_{\mathbb{R}^{3}}|Q|_\infty|u_{\varepsilon,p}|^{p^{*}}dx\right)\\
	&=& \left(\frac{1}{p}-\frac{1}{p^{*}}\right)\frac{\left(S_{p}^{\frac{3}{p}}+O\left(\varepsilon^{\frac{3-p}{p-1}}\right)\right)^{\frac{p^{*}}{p^{*}-p}}}{\left[\left(S_{p}^{\frac{3}{p}}+O\left(\varepsilon^{\frac{3}{p-1}}\right)\right)|Q|_{\infty}\right]^{\frac{p}{p^{*}-p}}} \\
	& =& \frac{1}{3}S_{p}^{\frac{3}{p}}|Q|_{\infty}^{-\frac{3-p}{p}}+O(\varepsilon^{\frac{3-p}{p-1}}).
\end{eqnarray*}According to Lemma \ref{lem0003.1},  there   exists  a unique $t_{\varepsilon} > 0$ such that $\widetilde{\ell}(t_{\varepsilon}) =\max\limits_{t>0} \mathcal{\widetilde J}(t u_{\varepsilon,p})>0$ and $\widetilde{\ell}^\prime(t_{\varepsilon})=0$. It is easy to verify that for some positive constants  $C_{1}, C_{2}$,
\begin{eqnarray}\label{CC}C_{1} \leq t_{\varepsilon} \leq C_{2}\quad \text{as $\varepsilon>0$ sufficiently small.}\end{eqnarray}
Let $\varepsilon$ be  sufficiently small and  $|x|\leq\varepsilon<\rho$.  Then, it follows  from \eqref{2.03} and \eqref{CC} that
\begin{eqnarray}\label{tU}
	t_{\varepsilon}u_{\varepsilon,p}=\frac{t_\varepsilon C\varepsilon^{\frac{3-p}{p(p-1)}}}{\left(\varepsilon^{\frac{p}{p-1}}+|x|^{\frac{p}{p-1}}\right)^{\frac{3-p}{p}}}\geq C\varepsilon^{\frac{p-3}{p}}.\end{eqnarray}
From \eqref{tU} and $(g_5$), we deduce that for any  $R>0$, there exists $s_0>0$ such that  $G(x,s)\geq R|s|^{\gamma}$ for $x\in\mathbb{R}^3$ and $|s|>s_0$. As a consequence,
\begin{eqnarray}\label{G>R}
	\label{2222}\int_{\mathbb{R}^3}G(x,{t}_{\varepsilon}u_{\varepsilon,p})dx\geq R{t}^\gamma_\varepsilon\int_{\mathbb{R}^3}|u_{\varepsilon,p}|^\gamma dx\geq RC\int_{\mathbb{R}^3}|u_{\varepsilon,p}|^\gamma dx.\end{eqnarray}
Then, by the assumption $(Q)$ and \eqref{G>R}, we derive that
\begin{eqnarray*}
	&& \max\limits_{t>0} \mathcal{\widetilde J}(tu_{\varepsilon,p})
	\\	&\leq& 		\!	\!\! \max_{t>0}\left(\frac{t^{p}}{p}\!\int_{\mathbb{R}^{3}}|\nabla u_{\varepsilon,p}|^{p}dx\!-	\!\frac{t^{p^{*}}}{p^{*}}\!\int_{\mathbb{R}^{3}}|Q|_\infty|u_{\varepsilon,p}|^{p^{*}}dx\right)\!+\!\frac{t_\varepsilon^{p^{*}}}{p^{*}}\!\!\int_{\mathbb{R}^{3}}(|Q|_\infty\!\!-\!Q(x))|u_{\varepsilon,p}|^{p^{*}}dx\\&&+\frac{(q-1)t_{\varepsilon }^{\frac{qm}{q-1}}}{qm}\!\int_{\mathbb{R}^{3}}\!K(x)\phi_{u_{\varepsilon,p}}|u_{\varepsilon,p}|^{m}dx \!+\!\frac{t_{\varepsilon }^{p}}{p}\int_{\mathbb{R}^{3}}\!V(x)|u_{\varepsilon,p}|^{p}dx\!-\! \int_{\mathbb{R}^{3}}\!G(x,t_{\varepsilon }u_{\varepsilon,p})dx
	\\&\leq&	\!	\!\frac{1}{3}S_{p}^{\frac{3}{p}}|Q|_{\infty}^{-\frac{3-p}{p}}+	M(\varepsilon),
\end{eqnarray*}
where
\begin{eqnarray*} \label{2.23}
	M(\varepsilon)\!\!\!&=&\!\!\! O(\varepsilon^{\frac{3-p}{p-1}})+\frac{t_\varepsilon^{p^*}C}{p^*}\!\int_{\mathbb{R}^{3}}|x|^\beta|u_{\varepsilon,p}|^{p^{*}}dx+\frac{(q-1)t_\varepsilon^{\frac{qm}{q-1}}}{qm}\!\int_{\mathbb{R}^{3}}K(x)\phi_{u_{\varepsilon,p}}|u_{\varepsilon,p}|^{m}dx 	\\&&\!\!\!+\frac{t_{\varepsilon }^{p}}{p}\int_{\mathbb{R}^{3}}V(x)|u_{\varepsilon,p}|^{p}dx- RC\int_{\mathbb{R}^3}|u_{\varepsilon,p}|^\gamma dx.
\end{eqnarray*}	
Next, we have only to prove that
\begin{eqnarray} \label{2.023}
 M(\varepsilon)<0 \quad\text{as   $\varepsilon>0$ small enough}.	\end{eqnarray}
It follows from  $(g_5$) that $\gamma>\frac{3(p-1)}{3-p}$. Then,  by  \eqref{2.08}, we obtain
\begin{eqnarray} \label{2.23a}
	M(\varepsilon)= O(\varepsilon^{\frac{3-p}{p-1}})+O\left(\varepsilon^{\frac{\gamma(3-p)}{p(p-1)}}\right)-RC\varepsilon^{3-\frac{\gamma(3-p)}{p}}+\widetilde{M}(\varepsilon), \end{eqnarray}	
where	
\begin{eqnarray*}
	\label{aa1}
		\left.
	\begin{array}{ll}
	\widetilde{M}(\varepsilon)\!=\!\!\displaystyle \frac{t_\varepsilon^{p^*}C}{p^*}\!\!\int_{\mathbb{R}^{3}}\!|x|^\beta|u_{\varepsilon,p}|^{p^{*}}dx\!+\!\frac{(q\!-\!1)t_\varepsilon^{\frac{qm}{q-1}}}{qm}\!\!\int_{\mathbb{R}^{3}}\!K(x)\phi_{u_{\varepsilon,p}}|u_{\varepsilon,p}|^{m}dx 	\!+\!\frac{t_{\varepsilon }^{p}}{p}\!\int_{\mathbb{R}^{3}}\!\!V(x)|u_{\varepsilon,p}|^{p}dx.\end{array}
\right.
\end{eqnarray*}		
 When $1<p\leq\sqrt{3}$ and $\beta\geq\frac{3}{p-1}$,  from \eqref{CC}, \eqref{2.07},  \eqref{2.07x}, \eqref{2.09},  one has
\begin{eqnarray}\label{a-1}
	\left.
	\begin{array}{ll}
	\widetilde{M}(\varepsilon)&\leq O\left(\varepsilon^\frac{3\!+\!\beta}{p}\right)\!+\!C\varepsilon^\frac{3\!+\!\beta}{p}|\ln\varepsilon|\!+\!O\left(\varepsilon^{\frac{3}{p\!-\!1}}\right)\!+\!	O\left(\varepsilon^{\frac{3-p}{p-1}}\right)\!+\!C\varepsilon^p\!+\!O\left(\varepsilon^p\right)\!+\!C\varepsilon^p|\ln\varepsilon|\\[3mm]&~~~~+\left\{\begin{array}{ll}O\big(\varepsilon^{\frac{mq(3-p)}{p(p-1)(q-1)}}\big)+C\varepsilon^{\frac{4pq-3qm+pqm-3p}{p(q-1)}},&\quad1<p<\frac{3qm+4q-3}{qm+4q-3},\\O\big(\varepsilon^{\frac{4q-3}{p(q-1)}}\big)+C\varepsilon^{\frac{4q-3}{p(q-1)}}\big|\ln\varepsilon\big|^{\frac{4q-3}{3(q-1)}},&\quad p=\frac{3qm+4q-3}{qm+4q-3},\\O\big(\varepsilon^{\frac{mq(3-p)}{p(q-1)(p-1)}}\big),&\quad\frac{3qm+4q-3}{qm+4q-3}<p<3.\end{array}\right.
	\end{array}
\right.\label{2}
\end{eqnarray}
Using  $1<p\leq\sqrt{3},\frac{3p}{4p-3}<q<3,m>\frac{(4q-3)p}{3q},\frac{qm}{q-1}\leq\gamma<p^*,\beta\geq\frac{3}{p-1}$, by elementary computations, it can be checked that
\begin{equation}\label{a-11}
	\resizebox{0.9\hsize}{!}{$	\displaystyle	0\!<	\displaystyle3-\frac{\gamma(3-p)}{p} \leq	\displaystyle\min\left\{\frac{3\!+\!\beta}{p},\frac{3}{p\!-\!1},\frac{3\!-\!p}{p\!-\!1},p,\frac{mq(3\!-p)}{p(p\!-1)(q\!-1)},\frac{4q\!-\!3}{p(q\!-\!1)},\frac{4pq\!-3qm\!+pqm\!-\!3p}{p(q\!-1)}\right\}.$}
\end{equation}
By \eqref{2.23a}, \eqref{a-1} and \eqref{a-11}, we derive that \eqref{2.023} holds.
When $1<p\leq\sqrt{3}$ and $\beta<\frac{3}{p-1}$, from \eqref{CC}, \eqref{2.07},  \eqref{2.07x}, \eqref{2.09},  we deduce that
\begin{eqnarray} \label{a-2}	\left.
	\begin{array}{ll}
	\widetilde{M}(\varepsilon)&\leq O\left(\varepsilon^{\frac{3}{p-1}}\right)+C\varepsilon^\beta+	O\left(\varepsilon^{\frac{3-p}{p-1}}\right)+C\varepsilon^p+O\left(\varepsilon^p\right)+C\varepsilon^p|\ln\varepsilon|\\[3mm]&~~~~+\left\{\begin{array}{ll}O\big(\varepsilon^{\frac{mq(3-p)}{p(p-1)(q-1)}}\big)+C\varepsilon^{\frac{4pq-3qm+pqm-3p}{p(q-1)}},&\quad1<p<\frac{3qm+4q-3}{qm+4q-3},\\O\big(\varepsilon^{\frac{4q-3}{p(q-1)}}\big)+C\varepsilon^{\frac{4q-3}{p(q-1)}}\big|\ln\varepsilon\big|^{\frac{4q-3}{3(q-1)}},&\quad p=\frac{3qm+4q-3}{qm+4q-3},\\O\big(\varepsilon^{\frac{mq(3-p)}{p(q-1)(p-1)}}\big),&\quad\frac{3qm+4q-3}{qm+4q-3}<p<3.\end{array}\right.	\end{array}
\right.
\end{eqnarray}	
From $1<p\leq\sqrt{3},\frac{3p}{4p-3}<q<3,m>\frac{(4q-3)p}{3q},	\max\left\{\frac{(3-\beta )p}{3-p},\frac{qm}{q-1}\right\}\leq\gamma<p^*$,
by some elementary calculations, we have
\begin{equation}\label{a-22}
	\resizebox{0.9\hsize}{!}{$	\displaystyle	0\!<3-\!\frac{\gamma(3-p)}{p}\displaystyle\! \leq\!	\displaystyle\min\!\left\{\!\frac{3}{p\!-\!1},\beta,\frac{3\!-\!p}{p\!-\!1},p,\frac{mq(3\!-\!p)}{p(p\!-\!1)(q\!-\!1)},\frac{4q\!\!-\!3}{p(q\!\!-\!1)},\frac{4pq\!-\!3qm\!+\!pqm\!-\!3p}{p(q\!-\!1)}\!\right\}$}.
\end{equation}
By \eqref{2.23a}, \eqref{a-2} and \eqref{a-22},  \eqref{2.023} holds.
When $\sqrt{3}< p<3$ and $\beta\geq\frac{3}{p-1}$,  from \eqref{CC}, \eqref{2.07},  \eqref{2.07x}, \eqref{2.09}, one has
\begin{eqnarray} \label{a-3}		\left.
	\begin{array}{ll}\left.
	\begin{array}{ll}
	\widetilde{M}(\varepsilon)&\leq O\left(\varepsilon^\frac{3+\beta}{p}\right)+C\varepsilon^\frac{3+\beta}{p}|\ln\varepsilon|+O\left(\varepsilon^{\frac{3}{p-1}}\right)+	O\left(\varepsilon^{\frac{3-p}{p-1}}\right)\\[3mm]&~~~~+\left\{\begin{array}{ll}O\big(\varepsilon^{\frac{mq(3-p)}{p(p-1)(q-1)}}\big)+C\varepsilon^{\frac{4pq-3qm+pqm-3p}{p(q-1)}},&\quad1<p<\frac{3qm+4q-3}{qm+4q-3},\\O\big(\varepsilon^{\frac{4q-3}{p(q-1)}}\big)+C\varepsilon^{\frac{4q-3}{p(q-1)}}\big|\ln\varepsilon\big|^{\frac{4q-3}{3(q-1)}},&\quad p=\frac{3qm+4q-3}{qm+4q-3},\\O\big(\varepsilon^{\frac{mq(3-p)}{p(q-1)(p-1)}}\big),&\quad\frac{3qm+4q-3}{qm+4q-3}<p<3.\end{array}\right.	\end{array}
\right.	\end{array}
\right.
\end{eqnarray}	
Since  $\sqrt{3}< p<3,\frac{3p}{4p-3}<q<3,\beta\geq\frac{3}{p-1}$ and
\begin{equation*}
	\resizebox{0.74\hsize}{!}{$	\displaystyle\max\!\left\{	{p^*\!\!-\!\!\frac{mq}{(q-1)(p-1)}},\frac{p(4p-6)}{(3-p)(p-1)} ,\frac{3pq\!-\!3p\!-\!4q\!+\!3}{(3\!-\!p)(\!p\!-\!1)}, \frac{qm}{q-\!1}\right\}\!\leq\!\gamma\!<\!p^*$},
\end{equation*}
with a simple computation one finds that
\begin{equation}\label{a-33}
	\resizebox{0.9\hsize}{!}{$	\displaystyle0\!<\!3-\frac{\gamma(3-p)}{p}
		\displaystyle \!\leq\!\min\!\left\{\frac{3+\!\beta}{p},\!\frac{3}{p\!-1},\frac{3-
			\!p}{p-1},\!\frac{mq(3\!-\!p)}{p(p-\!1)(q\!-1)},\!\frac{4q\!\!-\!\!3}{p(q\!-\!1)},\!\frac{4pq\!-\!3qm\!+\!pqm\!-3p}{p(q-1)}\right\}.$}
\end{equation}
Using \eqref{2.23a}, \eqref{a-3} and \eqref{a-33}, we obtain that \eqref{2.023} holds.
When $\sqrt{3}<p<3$ and $\beta<\frac{3}{p-1}$,  from \eqref{CC}, \eqref{2.07},  \eqref{2.07x}, \eqref{2.09}, we derive
\begin{eqnarray} \label{a-4}	\left.
	\begin{array}{ll}
	\widetilde{M}(\varepsilon)&\leq O\left(\varepsilon^\frac{3}{p-1}\right)+C\varepsilon^{\beta}+	O\left(\varepsilon^{\frac{3-p}{p-1}}\right)\\[3mm]&~~~~+\left\{\begin{array}{ll}O\big(\varepsilon^{\frac{mq(3-p)}{p(p-1)(q-1)}}\big)+C\varepsilon^{\frac{4pq-3qm+pqm-3p}{p(q-1)}},&\quad1<p<\frac{3qm+4q-3}{qm+4q-3},\\O\big(\varepsilon^{\frac{4q-3}{p(q-1)}}\big)+C\varepsilon^{\frac{4q-3}{p(q-1)}}\big|\ln\varepsilon\big|^{\frac{4q-3}{3(q-1)}},&\quad p=\frac{3qm+4q-3}{qm+4q-3},\\O\big(\varepsilon^{\frac{mq(3-p)}{p(q-1)(p-1)}}\big),&\quad\frac{3qm+4q-3}{qm+4q-3}<p<3.\end{array}\right.	\end{array}
\right.
\end{eqnarray}	
From $\sqrt{3}<p<3,\frac{3p}{4p-3}<q<3$ and
\begin{equation*}
	\resizebox{0.8\hsize}{!}{$	\displaystyle\max\left\{	{p^*-\!\!\frac{mq}{(q\!-\!\!1)(p\!-1)}},\!\frac{p(4p\!-\!6)}{(3\!-\!p)(p\!-\!1)} ,\frac{3pq\!-\!3p\!-\!4q\!+\!3}{(3\!-\!p)(p\!-\!1)},\frac{(3-\beta )p}{3\!-\!p},\frac{qm}{q\!-1}\right\}\leq\gamma<p^*$},
\end{equation*}
a direct calculation shows that
\begin{equation}\label{a-44}
	\resizebox{0.9\hsize}{!}{$	\displaystyle	0<3-\frac{\gamma(3-p)}{p} \leq\displaystyle\min\left\{\frac{3}{p\!-\!1},\beta,\frac{3\!-\!p}{p\!-\!1},\frac{mq(3\!-\!p)}{p(p\!-\!1)(q\!-\!1)},\frac{4q\!-\!3}{p(q-1)},\frac{4pq\!-\!3qm\!+\!pqm\!-\!3p}{p(q-1)}\right\}.$}
\end{equation}
Applying \eqref{2.23a}, \eqref{a-4} and \eqref{a-44}, we derive that \eqref{2.023} holds. The proof is complete. $\square$

Let us give  the periodic system corresponding to the asymptotically periodic system \eqref{1002},
\begin{eqnarray} \label{Vp1002}
	\left\{
	\begin{array}{ll}
		-\Delta_pu+V_T(x)|u|^{p-2}u+K_T(x)\phi|u|^{m-2}u=g_T(x,u)+Q_T(x)|u|^{p^*-2}u&\quad\text{in}\:\mathbb{R}^3,\\-\Delta_q\phi=K_T(x)|u|^{m}&\quad\text{in}\:\mathbb{R}^3.
	\end{array}
	\right.
\end{eqnarray}
We denote by $\widetilde{\mathcal{N}}_{T}$ the Nehari manifold associated with system \eqref{Vp1002}, namely
\begin{eqnarray*}\widetilde{\mathcal{N}}_{T}=\{u\in E\setminus\{0\}:\langle  \widetilde{\mathcal{J}}_{T}'(u),u\rangle=0\},	\end{eqnarray*}
and set
\begin{eqnarray*}\widetilde{m}_T:=\inf\limits_{u\in \widetilde{\mathcal{N}}_{T}} \widetilde{\mathcal{J}}_T(u),\end{eqnarray*}
where
\begin{eqnarray*}\widetilde{\mathcal{J}}_T(u)&=&\frac{1}{p}\int_{\mathbb{R}^{3}}(|\nabla u|^p+V_T(x)|u|^p)dx+\frac{q-1}{qm}\int_{\mathbb{R}^{3}} K_T(x)\phi_{u}|u|^{m}dx-\int_{\mathbb{R}^{3}} G_T(x,u)dx\\&&-\frac{1}{p^*}\int_{\mathbb{R}^{3}}Q_T(x)|u|^{p^{*}}dx.
\end{eqnarray*}

\begin{rem}\label{rem 3.7}
 Arguing as in the proof of Remark $\ref{rem 3.5}$,  we can conclude that
	\begin{eqnarray*}\widetilde{m}\leq \inf\limits_{u\in E\backslash\{0\}} \max\limits_{t>0}\widetilde{\mathcal{J}}_T(tu)= \widetilde{m}_T. \end{eqnarray*}
\end{rem}
	\textbf{Proof of Theorem 1.3.} Lemma \ref{lemsub 3.3} implies that there exists a bounded sequence $\{u_n\}\subset E$ satisfying \eqref{c01}. Then  there exists $u\in  E$ such that, up to a subsequence,
	\begin{eqnarray*}
		\left.
		\begin{array}{ll}
			\begin{cases}u_n\rightharpoonup u&\text{in}~E,\\u_n\to u&\text{in}~L_{\text{loc}}^l(\mathbb{R}^3)\text{~for all~} l\in[p,p^*) ,\\u_n(x)\to u(x)&\text{a.e.}\:\text{in~}\mathbb{R}^3.\end{cases}
		\end{array}
		\right.
	\end{eqnarray*}
For every  $\varphi\in{C}_0^\infty(\mathbb{R}^3)$, by   Lemma \ref{lem2.0002}(iv), Lemma \ref{sjx2} and  \cite[Proposition 5.4.7]{MW2013}, we have  $\langle \mathcal{\widetilde J}^{\prime}(u),\varphi\rangle=0$.  It follows that $(u, {\phi}_u)$ is a solution of
	system \eqref{1002}. We will show that system \eqref{1002} has a  ground state solution.
	
	Assume $u=0$ and set
	\begin{eqnarray*}  \sup_{z\in\mathbb{R}^{3}}\int_{B_{\delta}(z)}|u_{n}|^{p}dx\to\widetilde{\mu}.\end{eqnarray*}
We first prove that $\widetilde{\mu}>0$. Assume for contradiction that $ \widetilde{\mu}=0$. Then, it follows from  \cite[Lemma 3.2]{DSW2023} that  $u_n\to0$ in $L^{l}(\mathbb{R}^{3})$ for all $l\in[p,p^*)$. Moreover, using   the  assumptions $(K)$, $(Q)$, $(V)$, Lemmas \ref{lem2.0002}, \ref{lems2.1} and \eqref{Sp}, we derive that	\begin{eqnarray*}&&\frac{|Q|_{\infty}}{S_{p}^{\frac{p*}{p}}}\left(\int_{\mathbb{R}^{3}}|\nabla u_{n}|^{p}dx\right)^{\frac{p*}{p}}\\
		&\geq&\!\!-|Q|_{\infty}\int_{\mathbb{R}^{3}}|u_n|^{p^{*}}dx\\
		&\geq&\!\!-\int_{\mathbb{R}^{3}}K(x)\phi_{u_{n}}|u_{n}|^{m}dx\!+\!\int_{\mathbb{R}^{3}}g(x,u_{n})u_{n}dx\!+\!\int_{\mathbb{R}^{3}}Q(x)|u_n|^{p^{*}}dx+o(1)\\
		&\geq&\!\!-\int_{\mathbb{R}^{3}}|\nabla u_{n}|^{p}dx+o(1).\end{eqnarray*}
This means that only the following two cases can occur: \begin{eqnarray*}\text{either~} \int_{\mathbb{R}^{3}}|\nabla u_{n}|^{p}dx+o(1)\geq S^{\frac{3}{p}}|Q|_{\infty}^{-\frac{3-p}{p}} \text{~or~} \int_{\mathbb{R}^{3}}|\nabla u_{n}|^{p}dx=o(1),\end{eqnarray*}
 then   by the  assumptions $(K),(V)$, Lemmas \ref{lem2.0002}, \ref{lems2.1}   and $u_n\to0$ in $L^{l}(\mathbb{R}^{3})$ for all $l\in[p,p^*)$,  we deduce from \eqref{c01} that
	\begin{eqnarray*}
	\widetilde{m}& =& \mathcal{\widetilde J}(u_{n})-\frac{1}{p^{*}}\langle \mathcal{\widetilde J}^{\prime}(u_{n}),u_{n}\rangle+o(1) \\
		& =& \left(\frac{1}{p}-\frac{1}{p^{*}}\right)\|u_{n}\|^{p}+\left(\frac{q-1}{qm}-\frac{1}{p^{*}}\right)\int_{\mathbb{R}^{3}} K(x)\phi_{u_{n}}|u_{n}|^{m}dx\\&& +\int_{\mathbb{R}^{3}}\left(\frac{1}{p^{*}}g(x,u_{n})u_{n}-G(x,u_{n})\right)dx+o(1) \\
		&  \geq& \left(\frac{1}{p}-\frac{1}{p^{*}}\right)\|u_{n}\|^{p}+o(1) \\
		& \geq& \frac{1}{3}S_p^{\frac{3}{p}}|Q|_{\infty}^{-\frac{3-p}{p}},
	\end{eqnarray*}
which contradicts to Lemma  \ref{lemsub3.6}.
If $\int_{\mathbb{R}^{3}}|\nabla u_{n}|^{p}dx=o(1)$, then sobolev inequality implies that $\|u_{n}\|=o(1)$. From \eqref{c01}, it is easy to see that  $\widetilde{m}= 0$, which yields a contradiction because $\widetilde{m}>0$.
Thus, we conclude that $\widetilde{\mu}>0.$  Then, up to a subsequence, there exists  a sequence $\{z_n\}\subset\mathbb{Z}^3$ such that
	\begin{eqnarray}\label{1c}
		\int_{B_\delta(z_n)}|u_n|^pdx\geq\frac{\widetilde{\mu}}{2}>0.\end{eqnarray}
We set $\widetilde{u}_n(x):= u_n(x+z_n)$. Then exists $\widetilde{u}\in E$ such that up to subsequence,
	\begin{eqnarray} \label{u}
	\left.
	\begin{array}{ll}
		\begin{cases}\widetilde{u}_n\rightharpoonup \widetilde{u}&\text{in}~E,\\\widetilde{u}_n\to \widetilde{u}&\text{in}~L_{\text{loc}}^l(\mathbb{R}^3)\text{~for all~} l\in[p,p^*) ,\\\widetilde{u}_n(x)\to \widetilde{u}(x)&\text{a.e.}\:\text{in~}\mathbb{R}^3.\end{cases}
	\end{array}
	\right.
\end{eqnarray}
According to \eqref{1c} and \eqref{u}, it is easy to see that  $\widetilde{u}\neq0$ and
 $\{z_n\} $ is unbounded.
 Then, up to a subsequence, $|z_n|\to\infty.$
Moreover, for every  $\varphi\in C_0^{\infty}(\mathbb{R}^{3})$, using \eqref{u}, Lemmas \ref{lem2.0002} and \ref{sjx2}, by \cite[Proposition 5.4.7]{MW2013}, Lemmas \ref{lemc2.5},  \ref{wdb zm} and \ref{lem2.4},  we deduce that
		\begin{eqnarray*}
			\langle\mathcal{\widetilde J}^{\prime}_{T}(\widetilde{u}),\varphi\rangle	\!\!\!&=&\!\!\!\!\int_{\mathbb{R}^{3}}\left(|\nabla \widetilde{u}|^{p-2}\nabla \widetilde{u}\nabla\varphi+V_T(x)|\widetilde{u}|^{p-2}\widetilde{u}\varphi \right)dx+\int_{\mathbb{R}^{3}}K_T(x)\widetilde{\phi}_{\widetilde{u}}|\widetilde{u}|^{m-2}\widetilde{u}\varphi dx\\&&-\int_{\mathbb{R}^{3}}g_T(x,\widetilde{u})\varphi dx-\int_{\mathbb{R}^{3}}Q_T(x)|\widetilde{u}|^{p^{*}-2}\widetilde{u}\varphi dx	\\&=&\!\!\int_{\mathbb{R}^{3}}\!\!\left(|\nabla \widetilde{u}_n|^{p-2}\nabla \widetilde{u}_n\nabla\varphi\!+\!\!V_T(x)|\widetilde{u}_n|^{p-2}\widetilde{u}_n\varphi \right)dx\!+\!\!\int_{\mathbb{R}^{3}}\!K_T(x)\widetilde{\phi}_{\widetilde{u}_n}|\widetilde{u}_n|^{m-2}\widetilde{u}_n\varphi dx\\&&\!\!-\int_{\mathbb{R}^{3}}g_T(x,\widetilde{u}_n)\varphi dx\!-\!\int_{\mathbb{R}^{3}}Q_T(x)|\widetilde{u}_n|^{p^{*}-2}\widetilde{u}_n\varphi dx+o(1)\\	&=&\int_{\mathbb{R}^{3}}\left(|\nabla u_{n}|^{p-2}\nabla u_n\nabla\varphi(\cdot-z_{n})+V_T(x)|u_{n}|^{p-2} u_{n}\varphi(\cdot-z_{n})\right)dx\\&&+\int_{\mathbb{R}^{3}}K_T(x)\widetilde{\phi}_{u_n}|u_{n}|^{m-2}u_n\varphi(\cdot-z_{n})dx-\int_{\mathbb{R}^{3}}g_T(x,u_{n})\varphi(\cdot-z_{n})dx\\&&-\int_{\mathbb{R}^{3}}Q_T(x)|u_n|^{p^{*}-2}u_n\varphi(\cdot-z_{n}) dx+o(1)\\	&=&\int_{\mathbb{R}^{3}}\left(|\nabla u_{n}|^{p-2}\nabla u_n\nabla\varphi(\cdot-z_{n})+V(x)|u_{n}|^{p-2} u_{n}\varphi(\cdot-z_{n})\right)dx\\&&+\int_{\mathbb{R}^{3}}K(x)\phi_{u_n}|u_{n}|^{m-2}u_n\varphi(\cdot-z_{n})dx-\int_{\mathbb{R}^{3}}g(x,u_{n})\varphi(\cdot-z_{n})dx\\&&-\int_{\mathbb{R}^{3}}Q(x)|u_n|^{p^{*}-2}u_n\varphi(\cdot-z_{n}) dx+o(1)\\
		&=&\langle\mathcal{\widetilde J^{\prime}}(u_{n}),\varphi(\cdot-z_{n})\rangle+o(1)\\
		&=&0.
	\end{eqnarray*}
Thus, $(\widetilde{u}, \widetilde{\phi}_{\widetilde{u}})$ is the solution of system \eqref{Vp1002}. By the  assumption $(Q)$ and the  Fatou lemma, using Lemmas \ref{lemggGG}, \ref{lemc2.4} and \ref{lemmaV},   we deduce that
	\begin{eqnarray*}
	\widetilde{m}_{T}&\leq&\mathcal{\widetilde J}_T(\widetilde{u})\!-\!\frac{q-1}{qm}\langle\mathcal{\widetilde J}^{\prime}_T(\widetilde{u}),\widetilde{u}\rangle\\	&\leq&\left(\frac{1}{p}-\frac{q-1}{qm}\right)\int_{\mathbb{R}^{3}}\left(|\nabla \widetilde{u}|^{p}+V_T(x)|\widetilde{u}|^{p} \right) dx\!+\!\int_{\mathbb{R}^{3}}\left(\frac{q-1}{qm}g_T(x,\widetilde{u}) \widetilde{u}-G_T(x,\widetilde{u})\right)dx\\&&+\left(\frac{q-1}{qm}-\frac{1}{p^*}\right)\int_{\mathbb{R}^{3}}Q_T(x)|\widetilde{u}|^{p^{*}}dx\\
	&\leq&\left(\frac{1}{p}\!-\!\frac{q-1}{qm}\right)\int_{\mathbb{R}^{3}}\left(|\nabla \widetilde{u}_n|^{p}\!+\!V_T(x)|\widetilde{u}_n|^{p} \right) dx\!+\!\!\int_{\mathbb{R}^{3}}\left(\frac{q\!-\!1}{qm}g_T(x,\widetilde{u}_n) \widetilde{u}_n\!-\!G_T(x,\widetilde{u}_n)\right)dx\\&&+\left(\frac{q-1}{qm}-\frac{1}{p^*}\right)\int_{\mathbb{R}^{3}}Q_T(x)|\widetilde{u}_n|^{p^{*}}dx+o(1)\\
	&=&\left(\frac{1}{p}-\frac{q-1}{qm}\right)\int_{\mathbb{R}^{3}}\left(|\nabla u_n|^{p}\!+\!V_T(x)|u_n|^{p} \right) dx\!+\!\int_{\mathbb{R}^{3}}\left(\frac{q\!-\!1}{qm}g_T(x,u_{n}) u_{n}\!-\!G_T(x,u_{n})\right)dx\\&&+\left(\frac{q-1}{qm}-\frac{1}{p^*}\right)\int_{\mathbb{R}^{3}}Q_T(x)|u_n|^{p^{*}}dx+o(1)
	\\&\leq&\left(\frac{1}{p}\!-\!\frac{q\!-\!1}{qm}\right)\int_{\mathbb{R}^{3}}\left(|\nabla u_n|^{p}\!+\!V(x)|u_n|^{p} \right) dx+\int_{\mathbb{R}^{3}}\left(\frac{q-1}{qm}g(x,u_{n}) u_{n}\!-\!G(x,u_{n})\right)dx\\&&+\left(\frac{q-1}{qm}-\frac{1}{p^*}\right)\int_{\mathbb{R}^{3}}Q(x)|u_n|^{p^{*}}dx+o(1)\\
	&=&\mathcal{\widetilde J}(u_{n})-\frac{q-1}{qm}\langle\mathcal{\widetilde J^{\prime}}(u_{n}),u_{n}\rangle+o(1)\\&=& \widetilde{m}.
\end{eqnarray*}
Together with  Remark \ref{sremark3.2}, this implies that $\mathcal{\widetilde J}_T(\widetilde{u})=\widetilde{m}_T=\widetilde{m}>0.$ According to  Lemma \ref{lem0003.1}, there exists a unique $t_{\widetilde{u}}>0$ such that $t_{\widetilde{u}}\widetilde{u}\in\widetilde{\mathcal{N}}.$ Moreover, we get
 \begin{eqnarray}\widetilde{m}\leq \mathcal{\widetilde J}(t_{\widetilde{u}}\widetilde{u})\leq \mathcal{\widetilde J}_T(t_{\widetilde{u}}\widetilde{u})\leq \mathcal{\widetilde J}_T(\widetilde{u})=\widetilde{m}.\end{eqnarray}
Then $\mathcal{\widetilde J}(t_{\widetilde{u}}\widetilde{u})=\widetilde{m}$. It follows from Lemma \ref{lems4.4} that $(t_{\widetilde{u}}\widetilde{u},\phi_{t_{\widetilde{u}}\widetilde{u}})$ is the ground  state solution of system \eqref{1002}.

Assume $u\neq0.$ Since $\widetilde{u}\in\widetilde{\mathcal{N}}$, we get  \begin{eqnarray}\label{11} \mathcal{\widetilde J}(u)\geq \widetilde{m}.\end{eqnarray} Applying the assumptions $(V),(Q)$, the Fatou lemma and   Lemma \ref{lemggGG}, we deduce that	
\begin{eqnarray}
	\left.
	\begin{array}{ll}
		\mathcal{\widetilde J}(u)&=\displaystyle\mathcal{J}(u)-\frac{q-1}{qm}\langle \mathcal{J}^{\prime}(u),u\rangle\\[3mm]
		&=\displaystyle\left(\frac{1}{p}-\frac{q-1}{qm}\right)\int_{\mathbb{R}^{3}}\|u\|^{p}+\int_{\mathbb{R}^{3}}\left(\frac{q-1}{qm}g(x,u) u-G(x,u)\right)dx\\[3mm]
		&~~~~+\displaystyle\left(\frac{q-1}{qm}-\frac{1}{p^*}\right)\int_{\mathbb{R}^3}Q(x)|u|^{p^*}dx\\[3mm]
		&\leq\displaystyle\left(\frac{1}{p}-\frac{q-1}{qm}\right)\int_{\mathbb{R}^{3}}\|u_{n}\|^{p}+\int_{\mathbb{R}^{3}}\left(\frac{q-1}{qm}g(x,u_{n}) u_{n}-G(x,u_{n})\right)dx\\[3mm]
		&~~~~+\displaystyle\left(\frac{q-1}{qm}-\frac{1}{p^*}\right)\int_{\mathbb{R}^3}Q(x)|u_n|^{p^*}dx+o(1)\\[3mm]
		&=\displaystyle\mathcal{J}(u_{n})-\frac{q-1}{qm}\langle \mathcal{J}^{\prime}(u_{n}),u_{n}\rangle+o(1)\\[3mm]
		&= \widetilde{m}.	\end{array}
	\right.\label{D121}\end{eqnarray}
Then, from \eqref{11} and \eqref{D121}, we conclude that $\mathcal{\widetilde J}(u)=\widetilde{m}$.
	Thus $(u, {\phi}_u)$ is a  ground state solution of system \eqref{1002}.  The proof is complete. $\square$

\end{document}